\documentclass{amsart}

\usepackage{hyperref}
\usepackage{amsrefs}
\usepackage{comment}
\usepackage{float} 
\usepackage[dvipsnames]{xcolor}
\usepackage{tikz}
\usepackage{multirow}

\theoremstyle{definition}
\newcommand{\tri}{\mathbin{\overline{\nabla}}}

\newcommand{\fub}{ }

\newcommand{\m}{{\mathbf m}}
\newcommand{\C}{{\mathbf c}}
\newcommand{\T}{{\mathbf t}}
\newcommand{\Sb}{{\mathbf S}}

\newcommand{\Gb}{{\mathbf G}}
\newcommand{\Tb}{{\mathbf T}}

\newcommand{\Sm}{\mathcal{S}}
\newcommand{\Tm}{\mathcal{T}}  
\newcommand{\OB}{\big[}
\newcommand{\CB}{\big]}
\newcommand{\mC}{C_{\m}}
\newcommand{\mV}{V_{\m}}
\newcommand{\mE}{E_{\m}}
\newcommand{\mF}{F_{\m}}

\newcommand{\BigComment}[1]{}

\newcommand{\n}{{\bf n}}
\renewcommand{\k}{{\bf k}}

\providecommand{\U}[1]{\protect\rule{.1in}{.1in}}
 
\newtheorem{exercise}{Exercise}
\newtheorem{solution}{Solution}

\begin{document}
\title[Exercises for \textit{A Hyper-Catalan Series Solution} ]{Exercises for \textit{A Hyper-Catalan Series Solution to Polynomial Equations, and the Geode} }

\author{Dean Rubine}

\date{\today}
\begin{abstract}
We present exercises with solutions related to \textit{A Hyper-Catalan Series Solution to Polynomial Equations, and the Geode}.
\end{abstract}

\maketitle


The paper \textit{A Hyper-Catalan Series Solution to Polynomial Equations, and the Geode} by N.J. Wildberger and Dean Rubine~\cite{Wildberger2025}  is aimed at a relatively mathematically sophisticated audience, so skips some details that might not be immediately apparent to students.
We go through those details here as exercises with solutions. 
We'll refer to locations within the paper as WR followed by a section, page, theorem or equation.

\begin{exercise}  Solve $0=1-\alpha+t\alpha^2$ with the quadratic formula (WR eqn (2)).
\end{exercise}
\begin{solution}  Of course,
$$\alpha= \frac 1 {2t} (1 \pm \sqrt{1 - 4t }) \quad \checkmark$$
\end{solution}
\begin{exercise}  Show the series root $\alpha$ of $0=1-\alpha+t\alpha^2$ is the Catalan generating series by using the binomial theorem with $n=1/2$ to expand the minus sign square root.

This is a derivation referenced in the quadratic section (WR end of section 3); it isn't particularly germane to the current paper so may be skipped.
\end{exercise}
\begin{solution}  
Newton's innovation was to allow binomial coefficients with fractional tops.
The binomial expansion is going to give us $\binom{1/2}{k}$; let's work that out first.   For a fractional top, we need to write the binomial coefficient with the falling power in the numerator.
$$\binom{n}{k} = \dfrac{n(n-1)(n-2) \cdots(n-k+1)}{k!}
$$ 
$$  \binom{1/2}{k} = \dfrac{(1/2)(-1/2)(-3/2) \cdots(3/2-k)}{k!}   $$ 
That has $k$ fractions in the numerator; let's clear them with a factor of $2^k$.
$$  \binom{1/2}{k} = \dfrac{(1)(-1)(-3) \cdots(3-2k)}{2^k k!}  =  \dfrac{(-1)^{k-1}(1)(3)(5) \cdots(2k-3)}{2^k k!}   $$ 
where now we have the product of $k-1$ consecutive odd natural numbers; we'll compute that by dividing a factorial by a product of even numbers.
$$ (1)(3)(5) \cdots(2k-3) = \dfrac{ (2k-2)!}{(2)(4)(6) \dots (2k-2)} =  \dfrac{ (2k-2)!}{ 2^{k-1} (k-1)! }
$$
where we've factored $2^{k-1}$ out of the product of $k-1$ even numbers.  Substituting,
$$  \binom{1/2}{k} 
= \dfrac{(-1)^{k-1}(2k-2)! }{2^{2k-1} k! (k-1)!}   $$ 
With that in hand, let's give it a go, starting with the binomial expansion.
$$\alpha= \frac 1 {2t} (1 \! -  \! (1 \! - \! 4t )^{ \frac 1 2} ) = \frac 2 {4t} \! \left(1 \! - \! \sum_{k\ge 0} \binom{\frac 1 2}{k} (- 4t )^k  \right)= \frac 1 {2t} \! +\!  2 \sum_{k \ge 0} \binom{\frac 1 2}{k} (- 4t )^{k-1}    $$

We substitute $n=k-1$, or $k=n+1$.  The $n=-1$ term is $\binom{1/2}{0}(- \frac{1}{2t})$ which annihilates the initial $1/2t$.
$$\alpha=  \frac{1}{2t} + \sum_{n \ge -1} + \binom{1/2}{n+1} (-1)^n 2^{2n+1}t^n   =  \sum_{n \ge 0} \dfrac{(-1)^{n}(2n)! }{2^{2n+1} (n+1)! n! }   (-1)^{n} 2^{2n+1}t^n
$$
$$ \alpha =  \sum_{n \ge 0} \dfrac{(2n)!}{n! (n+1)!} t^n = \sum_{n \ge 0} C_n t^n \quad\checkmark
$$
\end{solution}
\begin{exercise} \label{exer:tri} For multisets of triagons $M_1,M_2$, show $\Psi(\tri (M_1,M_2)) = t \, \Psi(M_1) \Psi(M_2)$.  (WR p. 387)
\end{exercise}
\begin{solution}
We have
$ \tri(M_1, M_2) = \OB \ \tri(r_1, r_2) : \ r_1 \in M_1, r_2 \in M_2 \ \CB $
and 
$\psi(\tri (r,s)) = t \, \psi(r) \psi(s)$
and
$\Psi(M) \equiv \sum_{r \in M} \psi(r)$. 
\begin{align*}
\Psi &(\tri (M_1,M_2))  
= \Psi( [ \ \tri(r_1, r_2) : \ r_1 \in M_1, r_2 \in M_2 \ ] ) 
= \! \! \sum_{r_1 \in M_1} \sum_{r_2 \in M_2} \psi(\tri(r_1,r_2))
\\ &
=  \sum_{r_1 \in M_1} \sum_{r_2 \in M_2} t \psi(r_1) \psi(r_2)
=  t \, \sum_{r_1 \in M_1}  \psi(r_1) \sum_{r_2 \in M_2} \psi(r_2)
=  t \, \Psi(M_1)  \Psi(M_2) \quad\checkmark
\end{align*}
\end{solution}
\begin{exercise} Apply $\Psi$ to both sides of the multiset equation  $\Tm =  \ \OB \  | \ \CB \ +  \tri(\Tm,\Tm)$ to show $\Tb=\Psi(\Tm)$ satisfies $\Tb=1+t\Tb^2$.$\ \ \ $  (WR p. 387)
\end{exercise}
\begin{solution} Using results stated  in the previous exercise,
    \begin{align*}
\Tm & =  \ \OB \  | \ \CB \ +  \tri(\Tm,\Tm) 
\\ \Psi(\Tm) & =  \Psi( \ \OB \  | \ \CB \ +  \tri(\Tm,\Tm) \ )
\\ \Psi(\Tm) & =  \Psi( \OB \  | \ \CB)  + \Psi( \tri(\Tm,\Tm) )
\\ \Psi(\Tm) & =  \psi(|)  + t \, \Psi( \Tm) \Psi(\Tm)
\\ \Tb & =  1  + t  \Tb^2 \quad\checkmark
\end{align*}
\end{solution}

\begin{exercise}  Derive WR Theorem 2 from WR Theorem 1.  That is, show that since $\alpha=\sum_{n \ge 0} C_n t^n $ is the zero of $1-\alpha+t\alpha^2$ that
$ x  =  \sum_{n \ge 0} C_n \dfrac{c_0^{n+1} c_2^n}{c_1^{2n+1}}$
is the zero of $f(x)=c_0 - c_1 x + c_2 x^2$.
\end{exercise} 
\begin{solution}
Let $t=c_0c_2/c_1^2$ so $ \displaystyle g(\alpha) = 1-\alpha+ \frac{c_0 c_2}{c_1^2} \alpha^2$. 
Note that
$$ c_0 \, g\left(\frac{c_1}{c_0} x\right)  = c_0 \left(  1-\frac{c_1}{c_0} x  + \frac{c_0 c_2}{c_1^2} \left(\frac{c_1}{c_0} x \right)^2 \right) = c_0 - c_1 x + c_2 x^2 = f(x)
$$
\end{solution}
We have $\alpha=\sum_{n \ge 0} C_n t^n $ as the zero of $g(\alpha)$, which gives the zero of $f(x)$ for $\alpha = \frac{c_1}{c_0} x$, so
 $$ 
x = \frac{c_0}{c_1} \alpha = \frac{c_0}{c_1} \sum_{n \ge 0} C_n \left(  \frac{c_0 c_2}{c_1^2}  \right) ^n = \sum_{n \ge 0} C_n  \frac{c_0^{n+1} c_2^n}{c^{2n+1} } \quad\checkmark
$$
\begin{exercise} Generalizing exercise \ref{exer:tri} to multisets of subdigons (WR p. 390), show:
    $$\Psi(\tri_k(M_1, M_2, \ldots, M_k) ) = t_k \,  \Psi(M_2) \Psi(M_3) \cdots \Psi(M_k)$$
\end{exercise}
\begin{solution}
Proceeding as before, we have 
$\tri_k(M_1, M_2, \ldots, M_k) = \OB \tri_k(s_1, s_2, \ldots, s_k) \! : \ s_1 \in M_1, s_2 \in M_2, \ldots, s_k \in M_k  \CB$
and
$\psi(\tri_k(s_1, s_2, \ldots, s_k))=t_k \, \psi(s_1) \psi(s_2) \cdots \psi(s_k)$ 
and as before
$\Psi(M) \equiv \sum_{r \in M} \psi(r)$. 
\begin{align*}
\Psi& (\tri_k(M_1, M_2, \ldots, M_k) ) 
\\& = \Psi( [ \tri_k(s_1, s_2, \ldots, s_k) \! : \ s_1 \in M_1, s_2 \in M_2, \ldots, s_k \in M_k  ] ) 
\\&
= \sum_{s_1 \in M_1} \sum_{s_2 \in M_2} \cdots \sum_{s_k \in M_k} \psi(\tri_k(s_1, s_2, \ldots, s_k) )
\\ &
= \sum_{s_1 \in M_1} \sum_{s_2 \in M_2} \cdots \sum_{s_k \in M_k} t_k \, \psi(s_1) \psi(s_2)\cdots \psi(s_k) 
\\&
= t_k \sum_{s_1 \in M_1} \psi(s_1) \sum_{s_2 \in M_2} \psi(s_2) \cdots \sum_{s_k \in M_k} \psi(s_k)
\\& 
=  t_k \, \Psi(M_1)  \Psi(M_2) \cdots \Psi(M_k) \quad\checkmark
\end{align*}
\end{solution}
\begin{exercise}
Show $\Psi$ applied to both sides of  $\Sm =  \ \OB \  | \ \CB \ +  \tri_2(\Sm,\Sm) +  \tri_3(\Sm,\Sm,\Sm)  +  \tri_4(\Sm,\Sm,\Sm,\Sm)  + \ldots $ yields
$ \Sb = 1 + t_2 \Sb ^2 + t_3 \Sb ^3  + t_4 \Sb^4 + \ldots $ for $\Sb=\Psi(\Sm)$.$ \ $ (WR eqn (7))

\end{exercise}
\begin{solution}
\begin{align*}
\Sm &=  \ \OB \  | \ \CB \ +  \tri_2(\Sm,\Sm) +  \tri_3(\Sm,\Sm,\Sm)  +  \tri_4(\Sm,\Sm,\Sm,\Sm)  + \ldots 
\\ 
\Psi(\Sm) &= \Psi( \ \OB \  | \ \CB \ +  \tri_2(\Sm,\Sm) +  \tri_3(\Sm,\Sm,\Sm)  +  \tri_4(\Sm,\Sm,\Sm,\Sm)  + \ldots ) 
\\&
=  \Psi( \ \OB \  | \ \CB) + \Psi( \tri_2(\Sm,\Sm) ) + \Psi(  \tri_3(\Sm,\Sm,\Sm) ) + \Psi( \tri_4(\Sm,\Sm,\Sm,\Sm))  + \ldots 
\\&
=  \psi(|) + t_2 \Psi(\Sm) \Psi(\Sm)  + t_3 \Psi(\Sm) \Psi(\Sm) \Psi(\Sm)  + t_4 (\Psi( \Sm))^4   + \ldots 
\\
\Sb &= 1 + t_2 \Sb^2 + t_3 \Sb^3 + t_4 \Sb^4 + \ldots \quad\checkmark
\end{align*}
\end{solution}
\begin{exercise}
Derive WR Theorem 4 from WR Theorem 3.  That is, given  $\displaystyle  \alpha = \Sb[t_2,t_3, \ldots] =\sum_{\m \ge 0} \mC \T^{\m} $ is the zero of 
$\displaystyle g(\alpha) = 1 - \alpha + \sum_{k\ge 2} t_k  \alpha ^k$, show that
 the series zero of
$\displaystyle f(x)=c_0 - c_1 x + \sum_{k\ge 2} c_k x^k$
 is 
$\displaystyle x =  \sum_{\m \ge 0} \ \mC \  \dfrac{ c_0^{\mV -1} } { c_1^{\mE} } \, \C^{\m}$.
\end{exercise}
\begin{solution}
Let $t_k = c_0^{k-1}c_k /c_1^k$.  As before,
\begin{align*}
c_0 \, g\left(\frac{c_1}{c_0} x \right) &= c_0 \left(1 - \left(\frac{c_1}{c_0} x\right)  + \sum_{k\ge 2}  \frac{c_0^{k-1}c_k}{c_1^k}  \left(\frac{c_1}{c_0} x\right) ^k\right)  
\\
 &= c_0 - c_1 x + \sum_{k\ge 2} c_k x ^k = f(x)
\end{align*}
We'll also note (see WR eqns (4), (5), (6)):
$$\T^{\m} = \prod_{k\ge 2} t_k^{m_k} =  \ \prod_{k\ge 2} \left( \frac{c_0^{k-1}c_k}{c_1^k} \right)^{m_k} = \frac{c_0^{\sum_{k\ge 2} (k-1)m_k } \prod_{k\ge 2}  c_k^{m_k}}{c_1^{\sum_{k\ge 2} km_k}}   = \frac{c_0^{V_{\m}-2}\C^{\m} }{c_1^{E_{\m}-1}}
$$
So, to solve $f(x)=0$, we solve $g(\alpha)=0$ where $\alpha=\frac{c_1}{c_0}x$ giving
$$ x= \frac{c_0}{c_1}\alpha = \frac{c_0}{c_1} \sum_{\m \ge 0} \mC \T^{\m} = \frac{c_0}{c_1} \sum_{\m \ge 0} \mC \frac{c_0^{V_{\m}-2}\C^{\m} }{c_1^{E_{\m}-1}}=  \sum_{\m \ge 0} \ \mC \  \dfrac{ c_0^{\mV -1} } { c_1^{\mE} } \, \C^{\m} \quad\checkmark
$$
\end{solution}

\begin{exercise}
    Apply Schuetz and Whieldon 2016~\cite{Schuetz2016} Lemma 3.2 
    to derive WR Theorem 5, the explicit expression of the hyper-Catalans:
    $$\mC =  \dfrac{( 2m_2 + 3m_3 + 4m_4 + \ldots )!}{(1 + m_2 + 2m_3 + 3m_4 + \ldots)!  m_2! \, m_3! \,  m_4! \cdots} 
    = \dfrac{(E_{\m}-1)!}{(V_{\m}-1)! \, \m !} $$
\end{exercise}
\begin{solution}
Schuetz and Whieldon's polynomials only have certain exponents $d_i$, so they
restrict their polygonal subdivisions to only $(d_i+1)$-gons.
For the general polynomial we
require the general polygon subdivision, so we set $d_j=j+1$,  and rename $k_j$ to $m_{j+1}$ to match our notation.  The sums are over all the $k_j$; we drop the upper bound $r$ according to our convention:
\begin{align*}
k &= \sum_{j=1}^r k_j = \sum_{j \ge 1} m_{j+1} = \sum_{k\ge 2} m_{k} = F_{\m}
\\ n &= \sum_{j=1}^r (d_j - 1) k_j = \sum_{j 
\ge 1} j m_{j+1} = \sum_{k\ge 2} (k-1) m_{k} = V_{\m}-2 
\\ \textrm{ And since }  V_{\m} -{} & E_{\m} + F_{\m} = 1,
\\
n&+k  = F_{\m} +  V_{\m}-2 = E_{\m}-1 
\end{align*}
\begin{align*}
a_{\lambda} &= \dfrac{1}{n+1} \binom{n+k}{k}\binom{k}{k_1, k_2, \ldots}
=   \dfrac{1}{n+1} \dfrac{ (n+k)! }{k! n!} \dfrac{k!}{k_1! k_2 ! \cdots} 
 =    \dfrac{ (n+k)! }{ (n+1)! k_1! k_2 ! \cdots} 
 \\&
 = \dfrac{(E_{\m}-1)!}{(V_{\m}-1)! \, \m !} = \mC \quad\checkmark
\end{align*}    
\end{solution}
That covers the details in the derivations; we'll continue with some more exercises just for fun.

\begin{exercise} 
Let $ \displaystyle \beta = \Sb[ft_2, ft_3, ft_4,\ldots]= \sum_{\m \ge 0} \mC \C^{\m}$ where $c_n=ft_n$ for $n \ge 2$ so now the power of $f$ counts the number of faces (in the subdigons associated with its accounting monomial).
Show that $\beta$ satisfies $\beta - 1 = \sum_{k\ge 2} f t_k \beta^k$,
the general geometric polynomial with explicit faces.
 \end{exercise}
\begin{solution}
By WR Theorem 3, the zero of $\displaystyle 1 - \alpha + \sum_{k\ge 2} c_n \alpha^n$ is $$\alpha = \sum_{\m \ge 0} \mC \C^{\m} = \beta . \quad\checkmark$$
\end{solution}

\begin{exercise} 
Let $ \displaystyle \beta = \sum_{m_2+m_3 \le 4} C[m_2, m_3] (ft_2)^{m_2} (ft_3)^{m_3}$ where now the sum is only over terms with 4 or fewer faces.
When we evaluate $\beta-1$,  $ t_2 f \beta^2$ or $t_3 f \beta^3$ and collect powers of $f$ we'll get a polynomial in $t_2$ and $t_3$ multiplying each power $f^F$.  

Make a table with row for each face layer from $f^0$ to $f^4$,
with columns showing the associated polynomials for
$\beta-1$,  $ t_2 f \beta^2$ and  $t_3 f \beta^3$.
Verify at each face layer $F$, 
$[f^F] (\beta - 1) = [f^F] ( t_2 f \beta^2 + t_3 f \beta ^3  )$ where $[f^F]$ extracts the polynomial coefficient of $f^F$.

\end{exercise}
\begin{solution}
We can explicitly write out $\beta$; it has 15 terms. 
We do in the $\beta-1$ column.  
Squaring it and cubing $\beta$ gets a bit messy, but we only have to worry about terms that total to three faces or fewer, giving at most four faces after multiplication by $f$.
\begin{center}
\begin{tabular}{ | p{.45cm} | p{3.7cm}| p{3.2cm} | p{3.4cm} |} 
\hline
 $f^F$  & $\beta-1 \qquad\qquad = {}$  & $t_2 f \beta^2  \qquad + {}$  & $ t_3 f \beta^3$
\\  \hline
$f^0$ & 0 & 0 & 0
\\ $f^1 $&$ t_{2}+ t_{3} $ & $ t_{2}$ & $t_{3}$
\\ $f^2 $&$  2  t_{2}^{2}  + 5  t_{2} t_{3}+ 3t_{3}^{2}   $&$   2  t_{2}^{2}  +2  t_{2} t_{3} $&$  3 t_{2} t_{3} + 3 t_{3}^{2}   $
\\ $f^3 $&  \mbox{ $5 t_{2}^{3} + 21 t_{2}^{2} t_{3} $}  \mbox{$\ \ {}+ 28 t_{2} t_{3}^{2} + 12 t_{3}^{3} $} 
&  $ 5 t_{2}^{3} +  12 t_{2}^{2} t_{3} +7 t_{2} t_{3}^{2}  $&$   9 t_{2}^{2} t_{3}+21  t_{2} t_{3}^{2}  + 12  t_{3}^{3}  $ 
\\ $f^4 $ & \mbox{$14 t_{2}^{4} + 84 t_{2}^{3} t_{3} + 180 t_{2}^{2} t_{3}^{2} $}  \mbox{     $ {}+ 165 t_{2} t_{3}^{3} + 55 t_{3}^{4} $}
&$ 14 t_{2}^{4} + 56 t_{2}^{3} t_{3} + 72 t_{2}^{2} t_{3}^{2} + 30 t_{2} t_{3}^{3} $&$
28 t_{2}^{3} t_{3} + 108 t_{2}^{2} t_{3}^{2} + 135 t_{2} t_{3}^{3} + 55 t_{3}^{4} $ \\
\hline
\end{tabular}
\end{center}

We generated $\beta$ up to four faces, and shown agreement 
at the five face layers $F=0,1,2,3,4$.
Can you generalize and prove the result?
\end{solution}
\begin{exercise}
There are fourteen or fifteen sequences, slices of the hyper-Catalans, that we found in OEIS~\cite{OEIS}. Find some more.
\end{exercise}
\begin{solution}
    You're on your own here.
\end{solution}
\begin{exercise}
Gessel\cite{Gessel2016} gives the following formulation for Lagrange Inversion.  Given $R(\alpha)$, we can
find $\phi(f(x))$ where $f(x)$ is a series that satisfies  $f(x) = xR(f(x))$ and $\phi(x)$ is any Laurent series as follows:
$$
[x^n] \phi(f(x)) = \frac 1 n [\alpha^{n-1}] \phi'(\alpha) (R(\alpha))^n 
$$
Use Gessel's formulation to determine the series zero of the general geometric polynomial.
\end{exercise}
\begin{solution}
Let $a=f(1).$  Substituting, we find $a=R(a)$ so $a$ is the solution of $-\alpha + R(\alpha) = 0$.
We choose $R(\alpha)=1+t_2 \alpha^2 + t_3 \alpha^3 + \ldots$ $= \sum_{k\ge 1} t_k \alpha^k$ where we set $t_1=1/\alpha$ so the first term is $1$. 
$a$ is then a root of 
$$
0 = 1 - \alpha +t_2 \alpha^2 + t_3 \alpha^3 + \ldots
$$
the general geometric polynomial.
We want $(R(\alpha))^n$ so we write the multinomial theorem:
$$
\left( \sum_{m=1}^{M} s_m \right)^n = \sum_{\substack{k_i \ge 0 \\ \sum_i k_i =n }} \binom{n}{k_1, k_2, \ldots, k_M } \prod_{m=1}^M s_m ^{k_m} 
$$
The sum is over all the ways to make $M$ natural numbers $k_i$ add to $n$.  We'll drop $k_i \ge 0$ from the notation.
$$
(R(\alpha))^n 
= \left(\sum_{m=1}^M t_m \alpha^m \right)^n 
=\sum_{\substack{ \sum_i k_i =n }} \binom{n}{k_1, k_2, \ldots, k_M } \prod_{m=1}^M \left(  t_m \alpha^m \right)^{k_m}  
$$
$$
=\sum_{\substack{ \sum_i k_i =n }} \binom{n}{k_1, k_2, \ldots, k_M } \left( \prod_{m=1}^M   t_m ^{k_m} \right) \alpha^{\sum_{m=1}^M m  \, k_m  }
$$
We can just choose $\phi(\alpha)=\alpha$ so $\phi'(\alpha)=1$.
Then we're after 
$$
 \frac 1 n [\alpha^{n-1}] (R(\alpha))^n
 =   \frac 1 n [\alpha^{n-1}] \sum_{\substack{\sum_i k_i =n }} \binom{n}{k_1, k_2, \ldots, k_M } \left( \prod_{m=1}^M   t_m ^{k_m}  \right) \alpha^{\sum_{m=1}^M m  \, k_m  }
$$
We have to deal with $t_1=1/\alpha$, which we see has the net effect of starting the $m$ sum and product both at $m=2$.
$$
[x^n]f(x) =   \frac 1 n [\alpha^{n-1}] \sum_{\substack{ \sum_i k_i =n }} \binom{n}{k_1, k_2, \ldots, k_M } \left( \prod_{m=2}^M   t_m ^{k_m} \right) \alpha^{\sum_{m=2}^M m  \, k_m  }
$$
$$
[x^n]f(x) =   \sum_{\substack{\sum_{m=2}^M m k_m =n-1\\ \sum_{i=1}^M k_i =n }}  \frac 1 n \binom{n}{k_1, k_2, \ldots, k_M } \prod_{m=2}^M   t_m ^{k_m} 
$$
This is starting to look like our answer.
We still have a $k_1$ lurking around;
that's going to correspond to vertices.

Remember we want $a=f(1)$ so we set $x=1$ which has the net effect that $a$ is the sum of the coefficients of $f$,
$$
a =f(1) = \sum_{n \ge 0} \ \ \ \sum_{\substack{\sum_{m=2}^M m k_m =n-1\\ \sum_{m=1}^M k_m =n }}  \frac 1 n \binom{n}{k_1, k_2, \ldots, k_M } \prod_{m=2}^M   t_m ^{k_m} 
$$
OK, let's sort out $k_1$ and $n$.  Converting to our notation,
$$\m \equiv  [m_2, m_3, m_4, \ldots] = [k_2, k_3, k_4, \ldots] $$
Recall
$\mV 
= 2+\sum_{k \ge 2}^{\fub}  (k-1)m_k$,
$\mE= 1 + \sum_{k\ge 2}^{\fub} k \, m_k$,
$\mF = \sum_{k\ge 2}^{\fub} m_k$.
From our sum $a = \ldots$  we have
$$
\sum_{m=2}^M m k_m =n-1, \qquad \sum_{m=1}^M k_m =n 
$$
The first equation is $n=\mE$.
Clearly $\mF=\sum_{m=2}^M k_m$ so the second equation is $k_1 + \mF = n= \mE$.
That makes $k_1=\mE-\mF=\mV-1$ by the Euler Characteristic  formula $V-E+F=1$.
Back to the solution, which now looks like
$$
a  = \sum_{\mE \ge 0} \ \ \ \sum_{\mV-\mE+\mF=1}  \frac 1 {\mE} \binom{\mE}{\mV-1, m_2, m_3, \ldots } \prod_{i\ge 2}  t_i ^{k_i} 
$$
That's all looking very familiar.  The combination of the two sums is, we know, a sum over natural vectors; we have
$$
a  = \sum_{\m \ge 0} \dfrac{(\mE-1)!}{(\mV-1) ! \,  m_2 ! \,  m_3 ! \,  \ldots } \T^{\m}  = \sum_{\m \ge 0} \mC \T^{\m}  
$$

That's our derivation from first principles of Lagrange Inversion.
\end{solution}

\begin{exercise}
In a followup paper~\cite{Rubine2025}, we uncover this recursion for the hyper-Catalans, where $\vec{j}=[0,0,\ldots,1]$ with $j-2$ zeros and the multinomial coefficient bottom $\k$ denotes the (finitely many) non-zero elements of $k_\n$, order irrelevant:
$$
C_{\m} = \sum_{j\ge 2} \ \ \sum_{\substack{  \sum_{\n \ge 0} k_{\n} \n = \m - \vec{j} \\ \sum_{\n \ge 0} k_{\n} =j}} \binom{j}{\bf k} \prod_{\n \ge 0}  C_{\n} ^{k_{\n}} 
$$
Start from $C[\,]=1$ and generate some elements of the hyper-Catalan array.
\end{exercise}
\begin{solution}
Let's try the next Catalan number, $C[1]$.
$$
C[1] = \sum_{j\ge 2} \ \ \sum_{\substack{  \sum_{\n \ge 0} k_{\n} \n = [1] - \vec{j} \\ \sum_{\n \ge 0} k_{\n} =j}} \binom{j}{\bf k} \prod_{\n \ge 0}  C_{\n} ^{k_{\n}} 
$$
There's only one $j$, namely $j=2$, that gives a valid natural vector for $[1]-\vec{j}$.
$$
C[1] = \sum_{\substack{  \sum_{\n \ge 0} k_{\n} \n = [\,]  \\ \sum_{\n \ge 0} k_{\n} =2}} \binom{2}{\bf k} \prod_{\n \ge 0}  C_{\n} ^{k_{\n}} 
$$
All the non-null $\n$ have at least one non-zero element, so can't participate in the first sum, which says $k_{\n}=0$ for $\n \not = [\,]$. So the sum is over a single $\n=[\,]$; we conclude $k[\,]=2$, extending the unbounded subscript notation to the $k$ vector.
$$ 
C[1] = \binom{2}{2} ( C[\,] ) ^2 = 1
$$
Great.  What about $C[2]$?
$$
C[2] = \sum_{j\ge 2} \ \ \sum_{\substack{  \sum_{\n \ge 0} k_{\n} \n = [2] - \vec{j} \\ \sum_{\n \ge 0} k_{\n} =j}} \binom{j}{\bf k} \prod_{\n \ge 0}  C_{\n} ^{k_{\n}} 
$$
As before only $j=2$ enters because $\m=[2]$ only has a non-zero $m_2$.
$$
C[2] = \sum_{\substack{  \sum_{\n \ge 0} k_{\n} \n = [1]  \\ \sum_{\n \ge 0} k_{\n} =2}} \binom{2}{\bf k} \prod_{\n \ge 0}  C_{\n} ^{k_{\n}} 
$$
Any type $\n$ with $n_2 > 1$ or $n_k>0$ for $k \ge 3$ can't satisfy the first sum, which must involve only $\n=[\,]$ and $\n=[1]$.  The null vector can't contribute to the first sum, so we conclude $k[1]=1$ so $k[\,]=1$ in the single term.
$$
C[2] = \binom{2}{1,1}  C[\,] ^{k[\,]}  C[1] ^{k[1]}  = 2 \cdot 1^1 \cdot 1^1 =2
$$
Cranking out those Catalans!  One more, $C[3]$.  Again $j=2$, and $[3]-\vec{j}=[2]$.  By similar reasoning we can only have linear combinations $k_{[1]} [1] + k_{[2]}[2]=[2]$ so either $k[1]=2, k[2]=0$ or $k[1]=1$ and $k[2]=1$.
In both cases we get $k[0]=0$.
We have a sum of two terms:
$$
C[3] = \binom{2}{0,2,0}  C[\,]^0 C[1]^2 C[2]^0  + \binom{2}{0,1,1}  C[\,] ^{0}  C[1] ^{1} C[2]^1 =1 + 2(2^1)=5
$$
All right, let's try a Fuss number~\cite{Fuss1795}, namely $C[0,1]$. 
$$
C[0,1] = \sum_{j\ge 2} \ \ \sum_{\substack{  \sum_{\n \ge 0} k_{\n} \n = [0,1] - \vec{j} \\ \sum_{\n \ge 0} k_{\n} =j}} \binom{j}{\bf k} \prod_{\n \ge 0}  C_{\n} ^{k_{\n}} 
$$
We see only $j=3$ yields a natural vector, $\m-\vec{j}=[\,]$.
$$
C[0,1] = \sum_{\substack{  \sum_{\n \ge 0} k_{\n} \n = [\,] \\ \sum_{\n \ge 0} k_{\n} =3}} \binom{3}{\bf k} \prod_{\n \ge 0}  C_{\n} ^{k_{\n}} 
$$
The first sum tells us only the null vector can have a non-zero weight; the second sum tells us that weight is $3$.
$$
C[0,1] = \binom{3}{3}  C[\,] ^{3}  = 1
$$
All right, next Fuss number, $C[0,2]$.
Again only $j=3$ works; we get $\m-\vec{j}=[0,1]$.
Only the null vector and $[0,1]$ can partipate in the linear combination. 
We get $k[0,1]=1$ so $k[\,]=2$ and we have
$$
C[0,2] = \binom{3}{2,1} C[\,]^2 C[0,1]^1 = 3
$$
Wow.  Let's dip in our toes a bit deeper with $C[1,1]$.
Now we have a sum over both $j=2$ and $j=3$.  For the first we're aiming for $[0,1]$; for the second $[1,0]$.
The $[0,1]$ gives a single term with $j=2, k[\,]=1, k[0,1]=1$ and the $[1,0]$ gives a single term with $j=3, k[\,]=2, k[1,0]=1$.  So
$$
C[1,1] =  \binom{2}{1,1} 1^1 1^1 + \binom{3}{2,1} 1^1 1^1 = 2+3=5  
$$

We could do some more, but so far this recursion works pretty well.
\end{solution}

\begin{exercise}
 Derive the Segner Catalan convolution recurrence $$C_{m+1}=\sum_{n=0}^m C_n C_{m-n}$$ from hyper-Catalan recurrence given in the previous exercise.
\end{exercise}
\begin{solution} We're given:
$$
C_{\m} = \sum_{j\ge 2} \ \ \sum_{\substack{  \sum_{\n \ge 0} k_{\n} \n = \m - \vec{j} \\ \sum_{\n \ge 0} k_{\n} =j}} \binom{j}{\bf k} \prod_{\n \ge 0}  C_{\n} ^{k_{\n}} 
$$
The Catalan case restricts the first sum to $j=2$, and the vectors $\m$ and $\n$ to a single element each we'll just write as $m$ and $n$ and revert to natural indices on $k$.
$$
C_m = \sum_{\substack{  \sum_{n \ge 0} n k_{n} = \ m -1, \\  \sum_{n \ge 0} k_{n} =2}} \binom{2}{k_0, k_1, k_2, \ldots} \prod_{n \ge 0}  C_{n} ^{k_{n}} 
$$
We either get a single non-zero $k_n=2$ or two non-zero $k_n=1$, so we can break the sum into two sums.
$$
C_m = \sum_{\substack{ 2n=m-1 }} \binom{2}{2}  C_{n} ^2 + 
\sum_{\substack{ n<p \\ n+p=m-1 }} \binom{2}{1,1}  C_{n}C_p 
$$
Let's substitute $m=m'+1$ then drop the primes:
$$
C_{m+1} = \sum_{\substack{ 2n=m }}   C_{n} ^2 + 
\sum_{\substack{ n<p \\ n+p=m }} 2 C_{n}C_p 
$$
That's pretty close to the way Segner had it in 1758. Let's change it to:
$$
C_{m+1} = \sum_{\substack{ n=p \\ n+p=m }}   C_{n}C_p + 
\sum_{\substack{ n<p \\ n+p=m }} C_{n}C_p + \sum_{\substack{ n>p \\ n+p=m }} C_{n}C_p 
$$
We get zero or one term in the first sum according to $m$ being odd or even.
Either way, we can combine the three sums into one:
$$
C_{m+1} = \sum_{\substack{n+p=m}} C_{n}C_p  
$$
or, in simpler notation:
$$
C_{m+1} = \sum_{n=0}^m  C_nC_{m-n} \quad\checkmark  
$$

\end{solution}
\begin{exercise}
Recreate the little Schroeder numbers, the Riordan numbers and the Cayley 1891 array with various instantiations of $\Sb$.    
\end{exercise}

\begin{solution}
We have $\Sb[t_2, t_3, t_4, \ldots] =\Sb\{t_k\} = \sum_{\m \ge 0} \mC t_2^{m_2} t_3^{m_3} t_4^{m_4} \cdots$; our job is to set values for the $t_k$ to create the various sequences.

The Little Schroeder numbers are the coefficients of $v$ in 
$\Sb[v,v^2,v^3,\ldots] =\Sb\{v^{k-1}\}$, A001003.
$$\Sb\{v^{k-1}\} = \ldots {}+ 20793 v^{8} + 4279 v^{7} + 903 v^{6} + 197 v^{5} + 45 v^{4} + 11 v^{3} + 3 v^{2} + v + 1
$$
The Riordan numbers A005043 are the coefficients of $e$ in:
$$\Sb\{e^{k}\} = \ldots {}+602 e^{10} + 232 e^{9} + 91 e^{8} + 36 e^{7} + 15 e^{6} + 6 e^{5} + 3 e^{4} + e^{3} + e^{2} + 1
$$
The Cayley numbers (A033282 dovetailed) are
\begin{align*}
&\Sb\{v^{k-1}f\} = \ldots  
\\& + v^{7} \left(429 f^{7} + 1287 f^{6} + 1485 f^{5} + 825 f^{4} + 225 f^{3} + 27 f^{2} + f\right)
\\& + v^{6} \left(132 f^{6} + 330 f^{5} + 300 f^{4} + 120 f^{3} + 20 f^{2} + f\right)
\\& + v^{5} \left(42 f^{5} + 84 f^{4} + 56 f^{3} + 14 f^{2} + f\right)
\\& + v^{4} \left(14 f^{4} + 21 f^{3} + 9 f^{2} + f\right)
\\& + v^{3} \left(5 f^{3} + 5 f^{2} + f\right)
+ v^{2} \left(2 f^{2} + f\right) 
+ vf
+ 1
\end{align*}
\end{solution}

I used Python, in particular sympy, to do these; there's some code in exercise \ref{ex:code} but lots more not shown.

\begin{exercise}

Do the same projections of $\Gb$, the reduced hyper-Catalan array, the Geode.
\end{exercise}
\begin{solution}
In sympy I formed the $t_k v^{k-1}e^k f$ form of $\Sb$ to 6 face levels and 8 degrees, which limits how deep I can go in these projections by setting the various variables to unity.  We look at a few terms of $\Gb=\Sb/\Sb_1$.

The Little Schroeder Geode polynomial is 
$$
\ldots {} + 706 v^{5} + 152 v^{4} + 34 v^{3} + 8 v^{2} + 2 v + 1
$$
The Riordan Geode polynomial is
$$
\ldots {} + 141 e^{7} + 55 e^{6} + 21 e^{5} + 9 e^{4} + 3 e^{3} + 2 e^{2} + 1
$$
The Cayley Geode polynomial is
\begin{align*}
...
& + v^{5} \left(132 f^{5} + 288 f^{4} + 216 f^{3} + 64  f^{2} + 6 f\right) 
\\& + v^{4} \left(42 f^{4} + 70  f^{3} + 35  f^{2} + 5  f\right) 
\\& + v^{3} \left(14 f^{3} + 16  f^{2} + 4 f\right) 
+ v^{2} \left(5  f^{2} + 3 f\right)
+ 2vf + 1
\end{align*}
\end{solution}
\begin{exercise}
The history alludes to a larger array, that counts binary trees of a given type (a degree sequence, $k_1$ unary notes, $k_2$ binary nodes, etc.).  It's attributed to Raney, 1960~\cite{Raney1960}, Tutte, 1964~\cite{Tutte1964}, and Kreweras, 1972~\cites{Kreweras1972,KrewerasEarnshaw2005}.
Let's call it the \textbf{Tutrank Array}, and denote it $T_{\k}$, where $\k=[k_1, k_2, k_3, \ldots]$, with a starting index of one.

We have $V_{\k} = 2 + k_2+2k_3 + 3k_4 + \ldots$, unchanged from the hyper-Catalan case.
With a $k_1$, we have $E_{\k} = 1 + k_1 + 2 k_2 + 3 k_3 + \ldots$ and $F_{\k}= k_1 + k_2 + k_3 + \ldots$.

The exercise is to learn about the Tutrank array.
Produce some tables analogous to the ones we've already seen to understand this array that contains the hyper-Catalans $C[m_2,m_3,\ldots]=T[0,m_2,m_3,\ldots]$.
\end{exercise}
\begin{solution}
We'll forget triagons and repurpose $\Tb[t_1,t_2,t_3,\ldots]$ as the generator of the Tutrank Array,
$$\Tb[t_1,t_2,t_3,\ldots]= \sum_{\k \ge 0} T_{\k} \T^{\k}$$ where in the $\k$ case we define $\T^{\k} \equiv t_1^{k_1} t_2^{k_2} \cdots$, indices starting at one.
$k_1$ counts two-gons in our polygon picture, 
 which we can draw as two edges between the same pair of vertices.
 A two-gon adds a face and an edge, but no vertices.
Adding two-gons, we can double up any edge, or triple up, or more, indefinitely, without increasing the number of vertices.
So we see immediately that vertex layering is
unbounded for $\Tb$, as face layering was in the 
hyper-Catalan case (and still is for $\Tb$).

As before, we bound face layerings by specifying the degree.
We can bound vertex layerings by specifying a maximum $k_1$.
I can't think of another way.

Let's call subdivided polygons that we can divide with two-gons as well as the rest of them \textbf{tubdigons}.  We have $T_{\k}$ as the number of tubdigons of type $\k$.

Let's debug the software with about the simplest thing we can do.   We'll recreate equation (9), which is the one line cubic formula.  It's up to three faces, with degree 3.  We should be able to generate the same things for the Tutrank, set $t_1=0$ and recover the  $\Sb_F[t_2,t_3]$ up to 3 faces.
OK; got that to work, let's look at some face layers.

\begin{align*}
\Tb[t_1ef,&t_2ve^2f,t_3v^3e^2f] \mod f^6  =
\\   f^{5} & (273 e^{15} t_{3}^{5} v^{10} + 1001 e^{14} t_{2} t_{3}^{4} v^{9} + 715 e^{13} t_{1} t_{3}^{4} v^{8} + 1430 e^{13} t_{2}^{2} t_{3}^{3} v^{8}  +
\\ &  1980 e^{12} t_{1} t_{2} t_{3}^{3} v^{7} + 990 e^{12} t_{2}^{3} t_{3}^{2} v^{7} + 660 e^{11} t_{1}^{2} t_{3}^{3} v^{6} +
\\ & 1980 e^{11} t_{1} t_{2}^{2} t_{3}^{2} v^{6} + 330 e^{11} t_{2}^{4} t_{3} v^{6} + \\ &  1260 e^{10} t_{1}^{2} t_{2} t_{3}^{2} v^{5} + 840 e^{10} t_{1} t_{2}^{3} t_{3} v^{5} + 42 e^{10} t_{2}^{5} v^{5} + 252 e^{9} t_{1}^{3} t_{3}^{2} v^{4} + \\ &  756 e^{9} t_{1}^{2} t_{2}^{2} t_{3} v^{4} + 126 e^{9} t_{1} t_{2}^{4} v^{4} + 280 e^{8} t_{1}^{3} t_{2} t_{3} v^{3} + 140 e^{8} t_{1}^{2} t_{2}^{3} v^{3} + \\ &  35 e^{7} t_{1}^{4} t_{3} v^{2} + 70 e^{7} t_{1}^{3} t_{2}^{2} v^{2} + 15 e^{6} t_{1}^{4} t_{2} v + e^{5} t_{1}^{5} ) 
+ {} 
\\ f^{4} & (55 e^{12} t_{3}^{4} v^{8} + 165 e^{11} t_{2} t_{3}^{3} v^{7} + 120 e^{10} t_{1} t_{3}^{3} v^{6} + 180 e^{10} t_{2}^{2} t_{3}^{2} v^{6} + \\ &  252 e^{9} t_{1} t_{2} t_{3}^{2} v^{5} + 84 e^{9} t_{2}^{3} t_{3} v^{5} + 84 e^{8} t_{1}^{2} t_{3}^{2} v^{4} + 168 e^{8} t_{1} t_{2}^{2} t_{3} v^{4} + \\ &  14 e^{8} t_{2}^{4} v^{4} + 105 e^{7} t_{1}^{2} t_{2} t_{3} v^{3} + 35 e^{7} t_{1} t_{2}^{3} v^{3} + \\ &  20 e^{6} t_{1}^{3} t_{3} v^{2} + 30 e^{6} t_{1}^{2} t_{2}^{2} v^{2} + 10 e^{5} t_{1}^{3} t_{2} v +  e^{4} t_{1}^{4}) +
{} \\ f^{3} & (12 e^{9} t_{3}^{3} v^{6} + 28 e^{8} t_{2} t_{3}^{2} v^{5} + 21 e^{7} t_{1} t_{3}^{2} v^{4} + \\ &  21 e^{7} t_{2}^{2} t_{3} v^{4} + 30 e^{6} t_{1} t_{2} t_{3} v^{3} + 5 e^{6} t_{2}^{3} v^{3} + 10 e^{5} t_{1}^{2} t_{3} v^{2} + \\ &  10 e^{5} t_{1} t_{2}^{2} v^{2} + 6 e^{4} t_{1}^{2} t_{2} v + e^{3} t_{1}^{3}) + 
{} \\ f^{2} & (3 e^{6} t_{3}^{2} v^{4} +   5 e^{5} t_{2} t_{3} v^{3} + 4 e^{4} t_{1} t_{3} v^{2} + 2 e^{4} t_{2}^{2} v^{2} + 3 e^{3} t_{1} t_{2} v + e^{2} t_{1}^{2}) +
{} \\  f & (e^{3} t_{3} v^{2} + e^{2} t_{2} v + e t_{1}) + 1 
\end{align*}
The natural layering for the Tutrank is edge layering, as (without further bounding) it produces polynomial (not power series) coefficients.
\begin{align*}
(-1 + \Tb&[t_1ef,t_2ve^2f,t_3v^2e^3f,\ldots]) \mod e^{10}  =
\\ &
e^{9} (f^{9} t_{1}^{9} + 36 f^{8} t_{1}^{7} t_{2} v + 84 f^{7} t_{1}^{6} t_{3} v^{2} + 252 f^{7} t_{1}^{5} t_{2}^{2} v^{2} + 126 f^{6} t_{1}^{5} t_{4} v^{3} \\& +  630 f^{6} t_{1}^{4} t_{2} t_{3} v^{3} + 420 f^{6} t_{1}^{3} t_{2}^{3} v^{3} + 126 f^{5} t_{1}^{4} t_{5} v^{4} + 504 f^{5} t_{1}^{3} t_{2} t_{4} v^{4} + 252 f^{5} t_{1}^{3} t_{3}^{2} v^{4} + 756 f^{5} t_{1}^{2} t_{2}^{2} t_{3} v^{4} \\& +  126 f^{5} t_{1} t_{2}^{4} v^{4} + 84 f^{4} t_{1}^{3} t_{6} v^{5} + 252 f^{4} t_{1}^{2} t_{2} t_{5} v^{5} + 252 f^{4} t_{1}^{2} t_{3} t_{4} v^{5} + 252 f^{4} t_{1} t_{2}^{2} t_{4} v^{5} + 252 f^{4} t_{1} t_{2} t_{3}^{2} v^{5} \\& +  84 f^{4} t_{2}^{3} t_{3} v^{5} + 36 f^{3} t_{1}^{2} t_{7} v^{6} + 72 f^{3} t_{1} t_{2} t_{6} v^{6} + 72 f^{3} t_{1} t_{3} t_{5} v^{6} + 36 f^{3} t_{1} t_{4}^{2} v^{6} + 36 f^{3} t_{2}^{2} t_{5} v^{6} \\& +  72 f^{3} t_{2} t_{3} t_{4} v^{6} + 12 f^{3} t_{3}^{3} v^{6} + 9 f^{2} t_{1} t_{8} v^{7} + 9 f^{2} t_{2} t_{7} v^{7} + 9 f^{2} t_{3} t_{6} v^{7} + 9 f^{2} t_{4} t_{5} v^{7} + f t_{9} v^{8}) + {}
\\& 
e^{8} (f^{8} t_{1}^{8} + 28 f^{7} t_{1}^{6} t_{2} v + 56 f^{6} t_{1}^{5} t_{3} v^{2} + 140 f^{6} t_{1}^{4} t_{2}^{2} v^{2} + 70 f^{5} t_{1}^{4} t_{4} v^{3} + 280 f^{5} t_{1}^{3} t_{2} t_{3} v^{3} \\& +  140 f^{5} t_{1}^{2} t_{2}^{3} v^{3} + 56 f^{4} t_{1}^{3} t_{5} v^{4} + 168 f^{4} t_{1}^{2} t_{2} t_{4} v^{4} + 84 f^{4} t_{1}^{2} t_{3}^{2} v^{4} + 168 f^{4} t_{1} t_{2}^{2} t_{3} v^{4} + 14 f^{4} t_{2}^{4} v^{4} \\& +  28 f^{3} t_{1}^{2} t_{6} v^{5} + 56 f^{3} t_{1} t_{2} t_{5} v^{5} + 56 f^{3} t_{1} t_{3} t_{4} v^{5} + 28 f^{3} t_{2}^{2} t_{4} v^{5} + 28 f^{3} t_{2} t_{3}^{2} v^{5} + 8 f^{2} t_{1} t_{7} v^{6} \\& +  8 f^{2} t_{2} t_{6} v^{6} + 8 f^{2} t_{3} t_{5} v^{6} + 4 f^{2} t_{4}^{2} v^{6} + f t_{8} v^{7}) +{}
\\& 
e^{7} (f^{7} t_{1}^{7} + 21 f^{6} t_{1}^{5} t_{2} v + 35 f^{5} t_{1}^{4} t_{3} v^{2} + 70 f^{5} t_{1}^{3} t_{2}^{2} v^{2} + 35 f^{4} t_{1}^{3} t_{4} v^{3} \\& +  105 f^{4} t_{1}^{2} t_{2} t_{3} v^{3} + 35 f^{4} t_{1} t_{2}^{3} v^{3} + 21 f^{3} t_{1}^{2} t_{5} v^{4} + 42 f^{3} t_{1} t_{2} t_{4} v^{4} + 21 f^{3} t_{1} t_{3}^{2} v^{4} + 21 f^{3} t_{2}^{2} t_{3} v^{4} \\& +  7 f^{2} t_{1} t_{6} v^{5} + 7 f^{2} t_{2} t_{5} v^{5} + 7 f^{2} t_{3} t_{4} v^{5} + f t_{7} v^{6}) +{}
\\& 
e^{6} (f^{6} t_{1}^{6} + 15 f^{5} t_{1}^{4} t_{2} v + 20 f^{4} t_{1}^{3} t_{3} v^{2} + 30 f^{4} t_{1}^{2} t_{2}^{2} v^{2} + 15 f^{3} t_{1}^{2} t_{4} v^{3} \\& +  30 f^{3} t_{1} t_{2} t_{3} v^{3} + 5 f^{3} t_{2}^{3} v^{3} + 6 f^{2} t_{1} t_{5} v^{4} + 6 f^{2} t_{2} t_{4} v^{4} + 3 f^{2} t_{3}^{2} v^{4} + f t_{6} v^{5}) +
\\&  
e^{5} (f^{5} t_{1}^{5} + 10 f^{4} t_{1}^{3} t_{2} v + 10 f^{3} t_{1}^{2} t_{3} v^{2} + 10 f^{3} t_{1} t_{2}^{2} v^{2} + 5 f^{2} t_{1} t_{4} v^{3} + 5 f^{2} t_{2} t_{3} v^{3} + f t_{5} v^{4}) + {}
\\& 
e^{4} (f^{4} t_{1}^{4} + 6 f^{3} t_{1}^{2} t_{2} v + 4 f^{2} t_{1} t_{3} v^{2} + 2 f^{2} t_{2}^{2} v^{2} + f t_{4} v^{3}) + {}
\\& 
e^{3} (f^{3} t_{1}^{3} + 3 f^{2} t_{1} t_{2} v + f t_{3} v^{2}) + {}
\\& 
e^{2} (f^{2} t_{1}^{2} + f t_{2} v) + {} \\ & e^1 ( f t_{1} )
\end{align*}
These are multiplied by $ev^2$ to get the actual Euler counts; e.g. the two-gon $eft_1$ has two edges, two vertices and one face.

Do we still get the Riordan numbers?
\begin{align*}
(-1 + \Tb&[e,e^2,e^3,\ldots]) \mod e^{12}  =
58786 e^{11} + 16796 e^{10} + 4862 e^{9} + 1430 e^{8} \\& + 429 e^{7} + 132 e^{6} + 42 e^{5} + 14 e^{4} + 5 e^{3} + 2 e^{2} + e
\end{align*}

Wow, we get the Catalans; that's a surprise.  If we did this right we could set $t_1=0$ and reevaluate to get the Riordan numbers.  We get
\begin{align*}
(-1 + \Tb&[0,e^2,e^3,\ldots]) \mod e^{12}  = 1585 e^{11} + 603 e^{10} + 232 e^{9} + 91 e^{8}\\& + 36 e^{7} + 15 e^{6} + 6 e^{5} + 3 e^{4} + e^{3} + e^{2}
\end{align*}
indeed the Riordan numbers.  Cool.

Wait.  Why are these the \textit{shifted} Riodan numbers?
Shouldn't we get the entire sequence with the initial `1 0'? 
Same goes for the Catalans for that matter.
Taking a look at my code; it's \begin{verbatim}
eS = expand( -1 + Tbetae(12))
\end{verbatim}
so this is really $\Tb-1$; I'll fix the above to reflect that.

What about the Tutrank Cayley array?  Doesn't really work because of the unboundedness.

\begin{align*}
    \Tb&[f,vf,v^2f,\ldots] =
 \\&{} \quad
 v^{6} (\ldots {}+ 226512 f^{7} + 60984 f^{6} + 13860 f^{5} + 2520 f^{4} + 336 f^{3} + 28 f^{2} + f) 
 \\&{}
 + v^{5} (\ldots {}+ 60984 f^{7} + 19404 f^{6} + 5292 f^{5} + 1176 f^{4} + 196 f^{3} + 21 f^{2} + f)
 \\&{}
 + v^{4} \ldots {}+ (13860 f^{7} + 5292 f^{6} + 1764 f^{5} + 490 f^{4} + 105 f^{3} + 15 f^{2} + f) 
 \\&{}
 + v^{3} (\ldots {}+ 2520 f^{7} + 1176 f^{6} + 490 f^{5} + 175 f^{4} + 50 f^{3} + 10 f^{2} + f)
 \\&{}
 + v^{2} (\ldots {}+ 336 f^{7} + 196 f^{6} + 105 f^{5} + 50 f^{4} + 20 f^{3} + 6 f^{2} + f)
 \\&{}
 + v (\ldots {}+ 28 f^{7} + 21 f^{6} + 15 f^{5} + 10 f^{4} + 6 f^{3} + 3 f^{2} + f)
\end{align*}
We get a power series for each coefficient.  Clearly we get triangular numbers for $v^1$ and pyramidal numbers for $v^2$.  $v^3$ and $v^4$ are in OEIS as well.

Let's generate some bounded vertex layers for completeness; I'll include the $1$ and the $ev^2$ factor for a change:
\begin{align*}
ev^2& \Tb[t_1ef,t_2ve^2f,t_3v^2e^3f,\ldots] \mod t_1^{6} \mod v^8 =
\\& 
v^{7} (126126 e^{16} f^{10} t_{1}^{5} t_{2}^{5} + 168168 e^{15} f^{9} t_{1}^{5} t_{2}^{3} t_{3} + 42042 e^{15} f^{9} t_{1}^{4} t_{2}^{5} \\&  + 36036 e^{14} f^{8} t_{1}^{5} t_{2}^{2} t_{4} + 36036 e^{14} f^{8} t_{1}^{5} t_{2} t_{3}^{2} + 60060 e^{14} f^{8} t_{1}^{4} t_{2}^{3} t_{3} \\& + 12012 e^{14} f^{8} t_{1}^{3} t_{2}^{5} + 5544 e^{13} f^{7} t_{1}^{5} t_{2} t_{5} + 5544 e^{13} f^{7} t_{1}^{5} t_{3} t_{4} \\& + 13860 e^{13} f^{7} t_{1}^{4} t_{2}^{2} t_{4} + 13860 e^{13} f^{7} t_{1}^{4} t_{2} t_{3}^{2} + 18480 e^{13} f^{7} t_{1}^{3} t_{2}^{3} t_{3} \\& + 2772 e^{13} f^{7} t_{1}^{2} t_{2}^{5} + 462 e^{12} f^{6} t_{1}^{5} t_{6} + 2310 e^{12} f^{6} t_{1}^{4} t_{2} t_{5} \\& + 2310 e^{12} f^{6} t_{1}^{4} t_{3} t_{4} + 4620 e^{12} f^{6} t_{1}^{3} t_{2}^{2} t_{4} + 4620 e^{12} f^{6} t_{1}^{3} t_{2} t_{3}^{2} \\& + 4620 e^{12} f^{6} t_{1}^{2} t_{2}^{3} t_{3} + 462 e^{12} f^{6} t_{1} t_{2}^{5} + 210 e^{11} f^{5} t_{1}^{4} t_{6} + 840 e^{11} f^{5} t_{1}^{3} t_{2} t_{5} \\& + 840 e^{11} f^{5} t_{1}^{3} t_{3} t_{4} + 1260 e^{11} f^{5} t_{1}^{2} t_{2}^{2} t_{4} + 1260 e^{11} f^{5} t_{1}^{2} t_{2} t_{3}^{2} + 840 e^{11} f^{5} t_{1} t_{2}^{3} t_{3} \\& + 42 e^{11} f^{5} t_{2}^{5} + 84 e^{10} f^{4} t_{1}^{3} t_{6} + 252 e^{10} f^{4} t_{1}^{2} t_{2} t_{5} + 252 e^{10} f^{4} t_{1}^{2} t_{3} t_{4} + 252 e^{10} f^{4} t_{1} t_{2}^{2} t_{4} \\& + 252 e^{10} f^{4} t_{1} t_{2} t_{3}^{2} + 84 e^{10} f^{4} t_{2}^{3} t_{3} + 28 e^{9} f^{3} t_{1}^{2} t_{6} + 56 e^{9} f^{3} t_{1} t_{2} t_{5}\\&  + 56 e^{9} f^{3} t_{1} t_{3} t_{4} + 28 e^{9} f^{3} t_{2}^{2} t_{4} + 28 e^{9} f^{3} t_{2} t_{3}^{2} + 7 e^{8} f^{2} t_{1} t_{6} + 7 e^{8} f^{2} t_{2} t_{5} \\& + 7 e^{8} f^{2} t_{3} t_{4} + e^{7} f t_{6}) +
{} \\ &
v^{6} (18018 e^{14} f^{9} t_{1}^{5} t_{2}^{4} + 16632 e^{13} f^{8} t_{1}^{5} t_{2}^{2} t_{3} + 6930 e^{13} f^{8} t_{1}^{4} t_{2}^{4} \\& + 2772 e^{12} f^{7} t_{1}^{5} t_{2} t_{4} + 1386 e^{12} f^{7} t_{1}^{5} t_{3}^{2} + 6930 e^{12} f^{7} t_{1}^{4} t_{2}^{2} t_{3} \\& + 2310 e^{12} f^{7} t_{1}^{3} t_{2}^{4} + 252 e^{11} f^{6} t_{1}^{5} t_{5} + 1260 e^{11} f^{6} t_{1}^{4} t_{2} t_{4} \\& + 630 e^{11} f^{6} t_{1}^{4} t_{3}^{2} + 2520 e^{11} f^{6} t_{1}^{3} t_{2}^{2} t_{3} + 630 e^{11} f^{6} t_{1}^{2} t_{2}^{4} \\& + 126 e^{10} f^{5} t_{1}^{4} t_{5} + 504 e^{10} f^{5} t_{1}^{3} t_{2} t_{4} + 252 e^{10} f^{5} t_{1}^{3} t_{3}^{2} \\& + 756 e^{10} f^{5} t_{1}^{2} t_{2}^{2} t_{3} + 126 e^{10} f^{5} t_{1} t_{2}^{4} + 56 e^{9} f^{4} t_{1}^{3} t_{5} \\& + 168 e^{9} f^{4} t_{1}^{2} t_{2} t_{4} + 84 e^{9} f^{4} t_{1}^{2} t_{3}^{2} + 168 e^{9} f^{4} t_{1} t_{2}^{2} t_{3}\\&  + 14 e^{9} f^{4} t_{2}^{4} + 21 e^{8} f^{3} t_{1}^{2} t_{5} + 42 e^{8} f^{3} t_{1} t_{2} t_{4}\\&  + 21 e^{8} f^{3} t_{1} t_{3}^{2} + 21 e^{8} f^{3} t_{2}^{2} t_{3} + 6 e^{7} f^{2} t_{1} t_{5} + 6 e^{7} f^{2} t_{2} t_{4} + 3 e^{7} f^{2} t_{3}^{2} + e^{6} f t_{5}) +
{} \\ &
 v^{5} (2310 e^{12} f^{8} t_{1}^{5} t_{2}^{3} + 1260 e^{11} f^{7} t_{1}^{5} t_{2} t_{3} + 1050 e^{11} f^{7} t_{1}^{4} t_{2}^{3} + 126 e^{10} f^{6} t_{1}^{5} t_{4} \\ & + 630 e^{10} f^{6} t_{1}^{4} t_{2} t_{3} +  420 e^{10} f^{6} t_{1}^{3} t_{2}^{3} + 70 e^{9} f^{5} t_{1}^{4} t_{4} + 280 e^{9} f^{5} t_{1}^{3} t_{2} t_{3} \\&  + 140 e^{9} f^{5} t_{1}^{2} t_{2}^{3} + 35 e^{8} f^{4} t_{1}^{3} t_{4}  + 105 e^{8} f^{4} t_{1}^{2} t_{2} t_{3} + 35 e^{8} f^{4} t_{1} t_{2}^{3} \\& + 15 e^{7} f^{3} t_{1}^{2} t_{4} + 30 e^{7} f^{3} t_{1} t_{2} t_{3} + 5 e^{7} f^{3} t_{2}^{3} + 5 e^{6} f^{2} t_{1} t_{4} + 5 e^{6} f^{2} t_{2} t_{3} + e^{5} f t_{4}) +
{} \\ &
 v^{4} (252 e^{10} f^{7} t_{1}^{5} t_{2}^{2} + 56 e^{9} f^{6} t_{1}^{5} t_{3} + 140 e^{9} f^{6} t_{1}^{4} t_{2}^{2} + 35 e^{8} f^{5} t_{1}^{4} t_{3} \\& + 70 e^{8} f^{5} t_{1}^{3} t_{2}^{2} + 20 e^{7} f^{4} t_{1}^{3} t_{3} \\& + 30 e^{7} f^{4} t_{1}^{2} t_{2}^{2} + 10 e^{6} f^{3} t_{1}^{2} t_{3} + 10 e^{6} f^{3} t_{1} t_{2}^{2} \\& + 4 e^{5} f^{2} t_{1} t_{3} + 2 e^{5} f^{2} t_{2}^{2} + e^{4} f t_{3}) + 
{} \\ &
v^{3} (21 e^{8} f^{6} t_{1}^{5} t_{2} + 15 e^{7} f^{5} t_{1}^{4} t_{2} + 10 e^{6} f^{4} t_{1}^{3} t_{2} + 6 e^{5} f^{3} t_{1}^{2} t_{2} + 3 e^{4} f^{2} t_{1} t_{2} + e^{3} f t_{2}) + 
{} \\ &
v^{2} (e^{6} f^{5} t_{1}^{5} + e^{5} f^{4} t_{1}^{4} + e^{4} f^{3} t_{1}^{3} + e^{3} f^{2} t_{1}^{2} + e^{2} f t_{1} + e)
\end{align*}
The $t_1$s clutter up many of the terms; we limited them to the fifth power here.
\end{solution}
\begin{exercise}
A very similar proof shows the Tutrank has a factorization that leads to the Jumbo Geode.
Generate some layerings for the Jumbo Geode.
\end{exercise}
\begin{solution}
\newcommand{\Jb}{{\bf J}}
The Jumbo Geode is
$$
\Jb[t_1,t_2,\ldots] = \frac{-1 +\sum_{\k \ge 0} \Tb_{\k} \T^{\k}  }{ \sum_{i\ge 1} t_i  }
$$
We'll bound our face layers to degree 4; we'll leave off $ev^2$.
\begin{align*}
\Jb[&ef, ve^2f, v^2e^3f, v^3e^4f] = 
\\& 
f^{5} (7084 e^{20} t_{4}^{5} v^{15} + 26565 e^{19} t_{3} t_{4}^{4} v^{14} + 19250 e^{18} t_{2} t_{4}^{4} v^{13} + 39270 e^{18} t_{3}^{2} t_{4}^{3} v^{13} \\&{}+ 13265 e^{17} t_{1} t_{4}^{4} v^{12} + 55930 e^{17} t_{2} t_{3} t_{4}^{3} v^{12} + 28560 e^{17} t_{3}^{3} t_{4}^{2} v^{12} + 37690 e^{16} t_{1} t_{3} t_{4}^{3} v^{11} + 19510 e^{16} t_{2}^{2} t_{4}^{3} v^{11} \\&{}+ 59840 e^{16} t_{2} t_{3}^{2} t_{4}^{2} v^{11} + 10200 e^{16} t_{3}^{4} t_{4} v^{11} + 25625 e^{15} t_{1} t_{2} t_{4}^{3} v^{10} + 39320 e^{15} t_{1} t_{3}^{2} t_{4}^{2} v^{10} + 40840 e^{15} t_{2}^{2} t_{3} t_{4}^{2} v^{10} \\&{}+ 27880 e^{15} t_{2} t_{3}^{3} t_{4} v^{10} + 1428 e^{15} t_{3}^{5} v^{10} + 8155 e^{14} t_{1}^{2} t_{4}^{3} v^{9} + 52115 e^{14} t_{1} t_{2} t_{3} t_{4}^{2} v^{9} + 17800 e^{14} t_{1} t_{3}^{3} t_{4} v^{9} \\&{}+ 9050 e^{14} t_{2}^{3} t_{4}^{2} v^{9} + 27840 e^{14} t_{2}^{2} t_{3}^{2} t_{4} v^{9} + 4760 e^{14} t_{2} t_{3}^{4} v^{9} + 16035 e^{13} t_{1}^{2} t_{3} t_{4}^{2} v^{8} + 16745 e^{13} t_{1} t_{2}^{2} t_{4}^{2} v^{8} \\&{}+ 34365 e^{13} t_{1} t_{2} t_{3}^{2} t_{4} v^{8} + 2940 e^{13} t_{1} t_{3}^{4} v^{8} + 11990 e^{13} t_{2}^{3} t_{3} t_{4} v^{8} + 6160 e^{13} t_{2}^{2} t_{3}^{3} v^{8} + 9900 e^{12} t_{1}^{2} t_{2} t_{4}^{2} v^{7} + 10165 e^{12} t_{1}^{2} t_{3}^{2} t_{4} v^{7} \\&{}+ 21340 e^{12} t_{1} t_{2}^{2} t_{3} t_{4} v^{7} + 7315 e^{12} t_{1} t_{2} t_{3}^{3} v^{7} + 1870 e^{12} t_{2}^{4} t_{4} v^{7} + 3850 e^{12} t_{2}^{3} t_{3}^{2} v^{7} + 1855 e^{11} t_{1}^{3} t_{4}^{2} v^{6} \\&{}+ 12040 e^{11} t_{1}^{2} t_{2} t_{3} t_{4} v^{6} + 2065 e^{11} t_{1}^{2} t_{3}^{3} v^{6} + 4235 e^{11} t_{1} t_{2}^{3} t_{4} v^{6} + 6545 e^{11} t_{1} t_{2}^{2} t_{3}^{2} v^{6} + 1155 e^{11} t_{2}^{4} t_{3} v^{6} + 2130 e^{10} t_{1}^{3} t_{3} t_{4} v^{5} \\&{}+ 3385 e^{10} t_{1}^{2} t_{2}^{2} t_{4} v^{5} + 3490 e^{10} t_{1}^{2} t_{2} t_{3}^{2} v^{5} + 2475 e^{10} t_{1} t_{2}^{3} t_{3} v^{5} + 132 e^{10} t_{2}^{5} v^{5} + 1115 e^{9} t_{1}^{3} t_{2} t_{4} v^{4} + 575 e^{9} t_{1}^{3} t_{3}^{2} v^{4} \\&{}+ 1845 e^{9} t_{1}^{2} t_{2}^{2} t_{3} v^{4} + 330 e^{9} t_{1} t_{2}^{4} v^{4} + 125 e^{8} t_{1}^{4} t_{4} v^{3} + 555 e^{8} t_{1}^{3} t_{2} t_{3} v^{3} + 300 e^{8} t_{1}^{2} t_{2}^{3} v^{3} + 55 e^{7} t_{1}^{4} t_{3} v^{2} \\&{}+ 120 e^{7} t_{1}^{3} t_{2}^{2} v^{2} + 20 e^{6} t_{1}^{4} t_{2} v + e^{5} t_{1}^{5}) + 
{} \\&  
f^{4} (969 e^{16} t_{4}^{4} v^{12} + 2907 e^{15} t_{3} t_{4}^{3} v^{11} + 2091 e^{14} t_{2} t_{4}^{3} v^{10} + 3213 e^{14} t_{3}^{2} t_{4}^{2} v^{10} \\&{}+ 1411 e^{13} t_{1} t_{4}^{3} v^{9} + 4522 e^{13} t_{2} t_{3} t_{4}^{2} v^{9} + 1547 e^{13} t_{3}^{3} t_{4} v^{9} + 2962 e^{12} t_{1} t_{3} t_{4}^{2} v^{8} + 1549 e^{12} t_{2}^{2} t_{4}^{2} v^{8} \\&{}+ 3185 e^{12} t_{2} t_{3}^{2} t_{4} v^{8} + 273 e^{12} t_{3}^{4} v^{8} + 1958 e^{11} t_{1} t_{2} t_{4}^{2} v^{7} + 2015 e^{11} t_{1} t_{3}^{2} t_{4} v^{7} + 2119 e^{11} t_{2}^{2} t_{3} t_{4} v^{7} \\&{}+ 728 e^{11} t_{2} t_{3}^{3} v^{7} + 591 e^{10} t_{1}^{2} t_{4}^{2} v^{6} + 2570 e^{10} t_{1} t_{2} t_{3} t_{4} v^{6} + 442 e^{10} t_{1} t_{3}^{3} v^{6} + 453 e^{10} t_{2}^{3} t_{4} v^{6} \\&{}+ 702 e^{10} t_{2}^{2} t_{3}^{2} v^{6} + 737 e^{9} t_{1}^{2} t_{3} t_{4} v^{5} + 783 e^{9} t_{1} t_{2}^{2} t_{4} v^{5} + 810 e^{9} t_{1} t_{2} t_{3}^{2} v^{5} + 288 e^{9} t_{2}^{3} t_{3} v^{5} + 421 e^{8} t_{1}^{2} t_{2} t_{4} v^{4} \\&{}+ 218 e^{8} t_{1}^{2} t_{3}^{2} v^{4} + 468 e^{8} t_{1} t_{2}^{2} t_{3} v^{4} + 42 e^{8} t_{2}^{4} v^{4} + 69 e^{7} t_{1}^{3} t_{4} v^{3} + 232 e^{7} t_{1}^{2} t_{2} t_{3} v^{3} + 84 e^{7} t_{1} t_{2}^{3} v^{3} \\&{}+ 34 e^{6} t_{1}^{3} t_{3} v^{2} + 56 e^{6} t_{1}^{2} t_{2}^{2} v^{2} + 14 e^{5} t_{1}^{3} t_{2} v + e^{4} t_{1}^{4})
+ {} \\& 
f^{3} (140 e^{12} t_{4}^{3} v^{9} + 315 e^{11} t_{3} t_{4}^{2} v^{8} + 224 e^{10} t_{2} t_{4}^{2} v^{7} + 231 e^{10} t_{3}^{2} t_{4} v^{7} \\&{}+ 146 e^{9} t_{1} t_{4}^{2} v^{6} + 319 e^{9} t_{2} t_{3} t_{4} v^{6} + 55 e^{9} t_{3}^{3} v^{6} + 199 e^{8} t_{1} t_{3} t_{4} v^{5} + 106 e^{8} t_{2}^{2} t_{4} v^{5} \\&{}+ 110 e^{8} t_{2} t_{3}^{2} v^{5} + 125 e^{7} t_{1} t_{2} t_{4} v^{4} + 65 e^{7} t_{1} t_{3}^{2} v^{4} + 70 e^{7} t_{2}^{2} t_{3} v^{4} + 34 e^{6} t_{1}^{2} t_{4} v^{3} + 77 e^{6} t_{1} t_{2} t_{3} v^{3} \\&{}+ 14 e^{6} t_{2}^{3} v^{3} + 19 e^{5} t_{1}^{2} t_{3} v^{2} + 21 e^{5} t_{1} t_{2}^{2} v^{2} + 9 e^{4} t_{1}^{2} t_{2} v + e^{3} t_{1}^{3}) +
{} \\&   
f^{2} (22 e^{8} t_{4}^{2} v^{6} + 33 e^{7} t_{3} t_{4} v^{5} + 23 e^{6} t_{2} t_{4} v^{4} + 12 e^{6} t_{3}^{2} v^{4} + 14 e^{5} t_{1} t_{4} v^{3} \\& {}+16 e^{5} t_{2} t_{3} v^{3} + 9 e^{4} t_{1} t_{3} v^{2} + 5 e^{4} t_{2}^{2} v^{2} + 5 e^{3} t_{1} t_{2} v + e^{2} t_{1}^{2}) + 
{} \\&  
f (4 e^{4} t_{4} v^{3} + 3 e^{3} t_{3} v^{2} + 2 e^{2} t_{2} v + e t_{1}) +
{} \\&   1
\end{align*}

Edge layers:
\begin{align*}
\Jb[&ef, ve^2f, v^2e^3f, v^3e^4f] =
\\ &1 
\\ {} & +
e^1 f t_{1}
\\ {}&  +
e^2 ( f^{2} t_{1}^{2} + 2 f t_{2} v )
\\ {}&  +
e^3 ( f^{3} t_{1}^{3} + 5 f^{2} t_{1} t_{2} v + 3 f t_{3} v^{2} )
\\ {}&  +
e^4 ( f^{4} t_{1}^{4} + 9 f^{3} t_{1}^{2} t_{2} v + 9 f^{2} t_{1} t_{3} v^{2} + 5 f^{2} t_{2}^{2} v^{2} + 4 f t_{4} v^{3} )
\\ {} & +
e^5( f^{5} t_{1}^{5} + 14 f^{4} t_{1}^{3} t_{2} v + 19 f^{3} t_{1}^{2} t_{3} v^{2} + 21 f^{3} t_{1} t_{2}^{2} v^{2} + 14 f^{2} t_{1} t_{4} v^{3} + 16 f^{2} t_{2} t_{3} v^{3} )
\\ {} & + {} \cdots 
\end{align*}

Jumbo Geode Riordan sequence:
$$
\Jb[e, e^2, e^3, \ldots] =  \ldots {} + 271e^6+ 85 e^{5} + 28 e^{4} + 9 e^{3} + 3 e^{2} + e + 1
$$
Doesn't seem to appear in OEIS.
\end{solution}

\begin{exercise}
    Let $\Tm$ be the unbounded multiset of tubidigons, layered by type $\k=[k_1, m_2, m_3, m_4, \ldots] \equiv [k_1;  \m]$:
    $$ \Tm = \sum _{\k \ge 0} \Tm_{\k}$$
    Generalizing the subdigon story, $\Tm$ satisfies:
    $$\Tm = [ \, | \, ] + \tri_1(\Tm) + \tri_2(\Tm, \Tm) + \tri_3(\Tm, \Tm, \Tm) +  {} \ldots$$
Generalizing and applying $\Psi$, $\Tb=\Psi(\Tm)$ and
$$\Tb = 1 + t_1 \Tb + t_2 \Tb^2 + t_3 \Tb^3 + {} \ldots$$
$$0 = 1 - (1 - t_1) \Tb + t_2 \Tb^2 + t_3 \Tb^3 + {} \ldots$$
Solve with the polynomial formula, equate the two forms for $\Tb$ and draw conclusions.
\end{exercise}
\begin{solution}
 Let $c_0=1$, $c_1=1-t_1$, $c_k=t_k$  so $\C=\T$. Then by the polynomial theorem,
 $$
 \Tb = \sum_{\m \ge 0} \mC \frac{1}{(1- t_1)^{\mE}} \T^{\m} 
$$
Ah, $1/(1-t_1)$ is the geometric progression.
$$
 \Tb = \sum_{\m \ge 0} \mC \left( \sum_{k \ge 0} t_1^ {k} \right) ^{\mE}  \T^{\m} 
$$
Writing the appropriate multinomial theorem,
$$
\Tb = \sum_{\m \ge 0} \mC \left( \sum_{\substack{j_i \ge 0 \\ \sum_i j_i = \mE }} \binom{\mE}{j_0, j_1, j_2, \ldots } \prod_{k  \ge 0} (t_1^k) ^{j_k} \right) \T^{\m} 
$$
$$
\Tb = \sum_{\m \ge 0} \mC \sum_{\substack{ \sum_i j_i = \mE }} \binom{\mE}{j_0, j_1, j_2, \ldots } t_1^{ \sum_{k  \ge 0} k j_k} \T^{\m} 
$$
$$T[k_1  ; \m] =  [t_1 ^{k_1} ; \T^{\m} ] \Tb $$
$$ T[k_1  ; \m] =\mC \sum_{\substack{ \sum_{i  \ge 0}  j_i = \mE  \\   \sum_{i  \ge 0} i j_i = k_1}} \binom{\mE}{j_0, j_1, j_2, \ldots } 
$$

$$ T[m_1  ; \m] =\mC \sum_{\substack{ \sum_{i  \ge 0}  j_i = \mE  \\   \sum_{i  \ge 0} i j_i = m_1}} \binom{\mE}{j_0, j_1, j_2, \ldots } 
$$
Changed notation a bit to free up $k_i$.

Recall when $\k=[k_1, k_2, k_3, \ldots]$, with a starting index of one, we have 
$$V_{\k} = 2 + k_2+2k_3 + 3k_4 + \ldots,$$
unchanged from the hyper-Catalan case, 
$$E_{\k} = 1 + k_1 + 2 k_2 + 3 k_3 + \ldots$$
$$F_{\k}= k_1 + k_2 + k_3 + \ldots$$

Gotta reconcile starting at $j_0$ not $j_1$.
Let $j_i=k_{i+1}$.
$$ T[m_1  ; \m] =\mC \sum_{\substack{ \sum_{i  \ge 0}  k_{i+1} = \mE  \\   \sum_{i  \ge 0} i k_{i+1} = m_1}} \binom{\mE}{k_1, k_2, \ldots } 
$$
$$
T[m_1  ; \m] =\mC \sum_{\substack{ \sum_{i  \ge 1}  k_{i} = \mE  \\   \sum_{i  \ge 1} (i-1) k_{i} = m_1}} \binom{\mE}{k_1, k_2, \ldots } 
$$
$$
T[m_1  ; \m] =\mC \sum_{\substack{ \sum_{i  \ge 1}  k_{i} = \mE  \\   \sum_{i  \ge 1} i k_{i} = m_1 + \mE }} \binom{\mE}{k_1, k_2, \ldots } 
$$
That's in the form where Fine's Identity applies,
$$
T[m_1  ; \m] =\binom{ \mE -1 + m_1  }{\mE-1 }  \mC = \binom{ \mE -1 + m_1  }{m_1} \mC 
$$

I thought that there was a a mistake and that the correct formula appeared to be $ T[m_1  ; \m] =\binom{ \mE -1 }{m_1} \mC $ but I was confusing $\mE$ and $E[m_1;\m]=m_1 + \mE$.

Did I do this the hard way?  Can we just factor out factors with $m_1$? Sure; piece of cake:
$$
C[m_2, m_3,\ldots] =  \dfrac{( 2m_2 + 3m_3 + 4m_4 + \ldots )!}{(1 + m_2 + 2m_3  + \ldots)!  m_2! \, m_3! \cdots} = \dfrac{( \mE-1)!}{ (\mV -1) ! \,  \m! \  } 
$$
$$T[m_1, m_2, m_3, \ldots] =  \dfrac{(m_1+ 2m_2 + 3m_3 + 4m_4 + \ldots )!}{(1 + m_2 + 2m_3 + 3m_4 + \ldots)! \, m_1! \, m_2! \, m_3! \,  m_4! \cdots} $$
$$ = \dfrac{(m_1 + \mE-1)!}{(\mV-1)! \, m_1! \, \m!} $$
$$ = \dfrac{(m_1 + \mE - 1)! } {  m_1!\,   (\mE-1) !} \cdot \dfrac{ (\mE-1) !} {(\mV-1)! \, \m!} $$
$$ T[m_1; \m] = \binom{m_1 + \mE - 1}{ m_1} \mC \quad\checkmark $$
Recapping what we did above, from:
$$\Tb = 1 + t_1 \Tb + t_2 \Tb^2 + t_3 \Tb^3 + {} \ldots$$
we get to:
$$
T[m_1  ; \m] =\mC \sum_{\substack{ \sum_{i  \ge 1}  k_{i} = \mE  \\   \sum_{i  \ge 1} i k_{i} = m_1 + \mE }} \binom{\mE}{k_1, k_2, \ldots } 
$$
If we didn't know Fine's Identity~\cite{Fine1988} equation (22.17), we could then conclude: 
$$
 \sum_{\substack{ \sum_{i  \ge 1}  k_{i} = \mE  \\   \sum_{i  \ge 1} i k_{i} = m_1 + \mE }} \binom{\mE}{k_1, k_2, \ldots } 
 =
\binom{m_1 + \mE - 1}{ m_1}
$$
where now $\mE$ and $m_1$ are just natural numbers.
Setting $n=m_1+\mE, k=\mE$, we've proven Fine's Identity using Wildberger's Polynomial Formula:
$$
 \sum_{\substack{ \sum_{i  \ge 1}  k_{i} = k  \\   \sum_{i  \ge 1} i k_{i} = n }} \binom{k}{k_1, k_2, \ldots } 
 =
\binom{n - 1}{k - 1}
$$
\end{solution}
\begin{exercise}

Try to understand NJW's conjecture about the Geode array.  Referring to the Bi-Tri Geode, he writes, ``that these numbers are counting ordered incomplete trees with $m_2$ binary nodes, $m_3$ ternary nodes and a single additional leaf node.
We conjecture that the natural generalization of this explains all the entries of the Geode array $G_{\m}$.''

\begin{figure}[H]
\centering
    \includegraphics[width=1\linewidth]{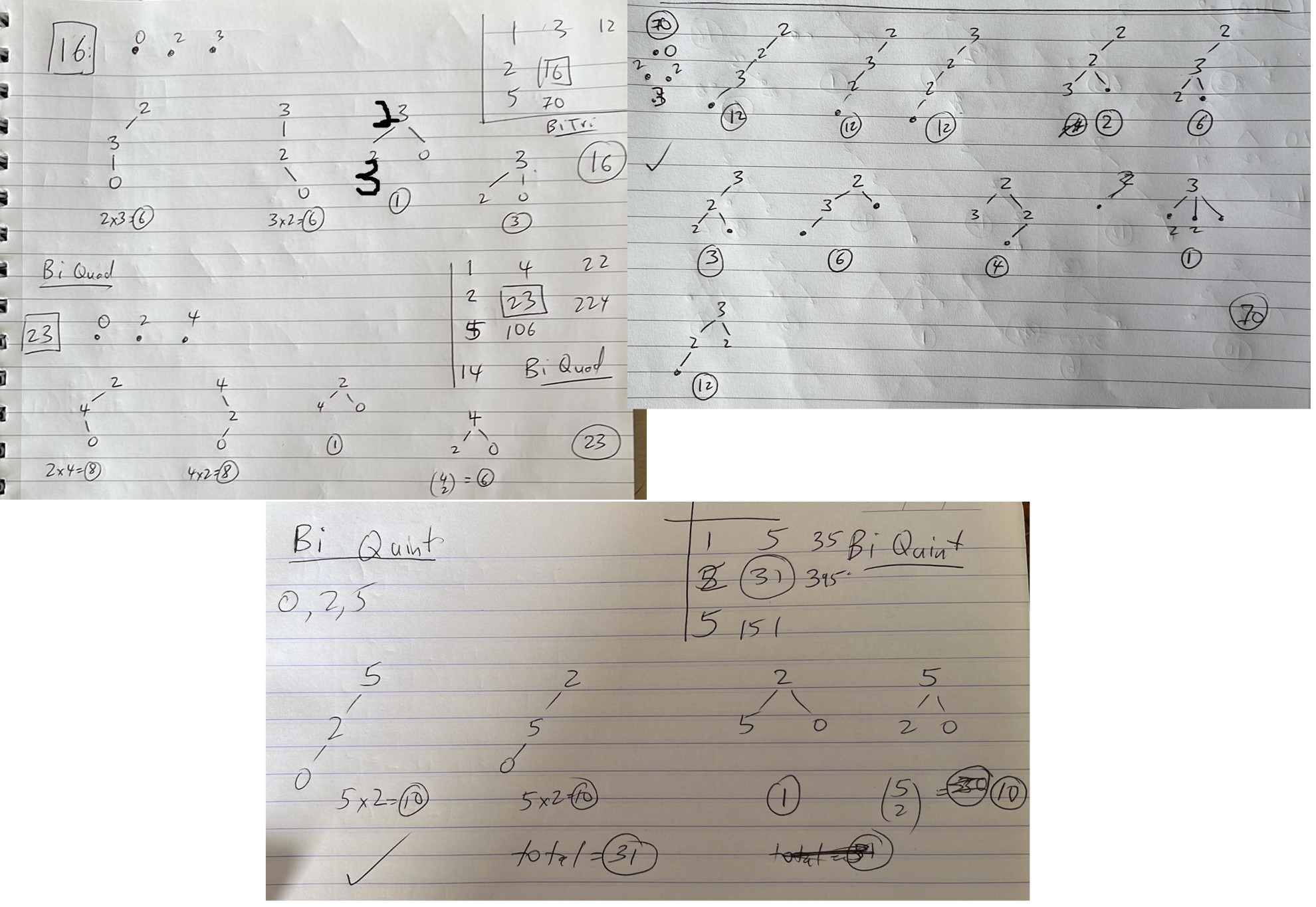} \caption{NJW's examples of counting Geode trees}
\label{fig:NG}
\end{figure}
\end{exercise}
\begin{solution}
Let's write the natural generalization as a strawman we can attack; we do that by adding an `etc.'

{\bf Conjecture G1. } $G_{\m}=G[m_2, m_3, m_4, \ldots]$ is the number of ordered incomplete trees with $m_2$ binary nodes, $m_3$ ternary nodes, etc., and a single additional leaf node.

Each of NJW's trees is schematic, generally representing multiple trees.
I suppose we should try the smallest first, NJW's 16, which appears to be $G[1,1]$.

I corrected the schema with a count of 1 to have the binary node as root; otherwise it seems this tree is counted in the 3 to the right of it.  Is this correction correct?

I wonder why we can't have an `ordered incomplete tree' where we put the 0 before the 2 or the 3, so  the trees:
\begin{figure}[H]
\centering
    \includegraphics[width=.5\linewidth]{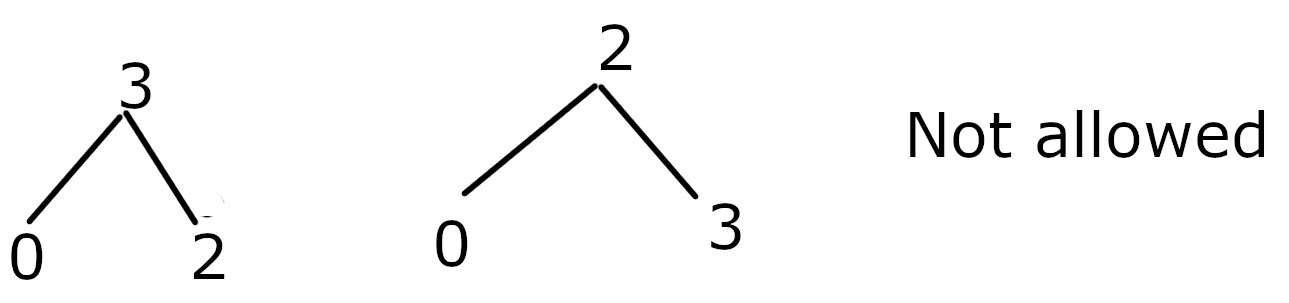} \caption{These trees don't seem to counted by the Geode array}
\label{fig:GD1}
\end{figure}

{\bf Conjecture G2. }  $G_{\m}$ is the number of ordered incomplete trees with $m_2$ binary nodes, $m_3$ ternary nodes, etc. and a single additional leaf node never to the left of any node with non-zero arity.

Let's go with it, tentatively. 

That gives me a different method than NJW for counting.
For each of the $\mC$ subdigons of type $\m$, figure out (count) which leaf nodes could be incomplete leaf nodes.
\begin{figure}[H]
\centering
    \includegraphics[width=1\linewidth]{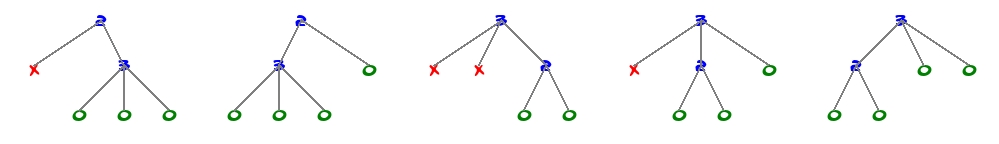} \caption{Dean's method of counting Geode array $C[1,1]=5, G[1,1]=16$}
\label{fig:DCG11}
\end{figure}
\begin{figure}[H]
\centering
    \includegraphics[width=1\linewidth]{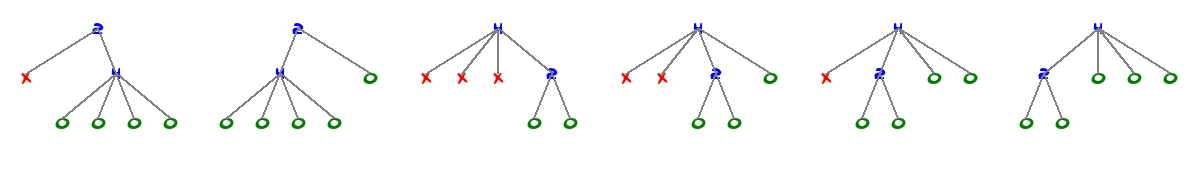} \caption{$C[1,0,1]=6, G[1,0,1]=23$}
\label{fig:DCG101}
\end{figure}
\begin{figure}[H]
\centering
    \includegraphics[width=1\linewidth]{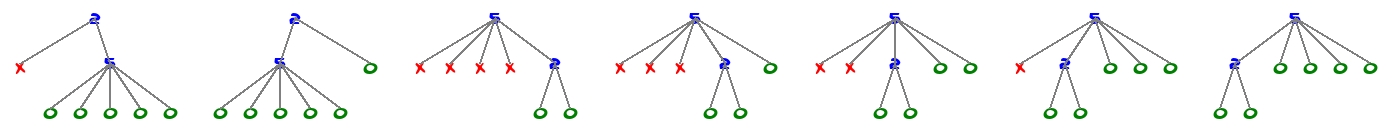} \caption{$C[1,0,0,1]=7, G[1,0,0,1]=31$}
\label{fig:DCG1001}
\end{figure}
\begin{figure}[H]
\centering
    \includegraphics[width=1\linewidth]{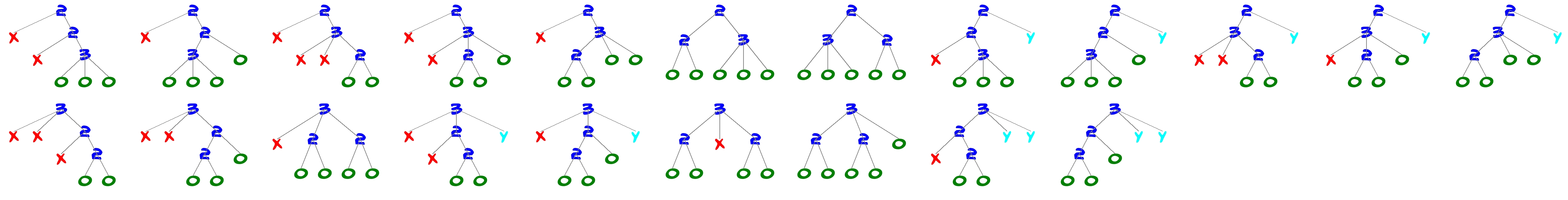} \caption{$C[2,1]=21, G[2,1]=70$}
\label{fig:DCG21}
\end{figure}

When I do this, applying the never left of a sibling rule to the additional leaf node, it works for $G[1,1]$, $G[1,0,1]$ and $G[1,0,0,1]$ but gives too many for $G[2,1]$.
When we don't count the the additional nodes at some nodes, indicated by cyan $Y$s, we get a correct count.  Those aren't random nodes; they correspond to configurations not counted by NJW in his enumeration. I'm not sure how to qualify the rule to exclude those nodes, let's try:

{\bf Conjecture G3. }   $G_{\m}$ is the number of ordered incomplete trees with $m_2$ binary nodes, $m_3$ ternary nodes, etc. and a single additional leaf node never to the left of any node with non-zero arity and never depth two or more less than maximum depth.  

At this point, NJW's description of a ``single additional leaf node'' appears to need serious modification. 

We should try another example.  They're getting out of hand calculation range.
Four faces seems hard; let's start from $C[1,2]=28$ and try to make $G[1,2]=110$.

\begin{figure}[H]
\centering
    \includegraphics[width=1\linewidth]{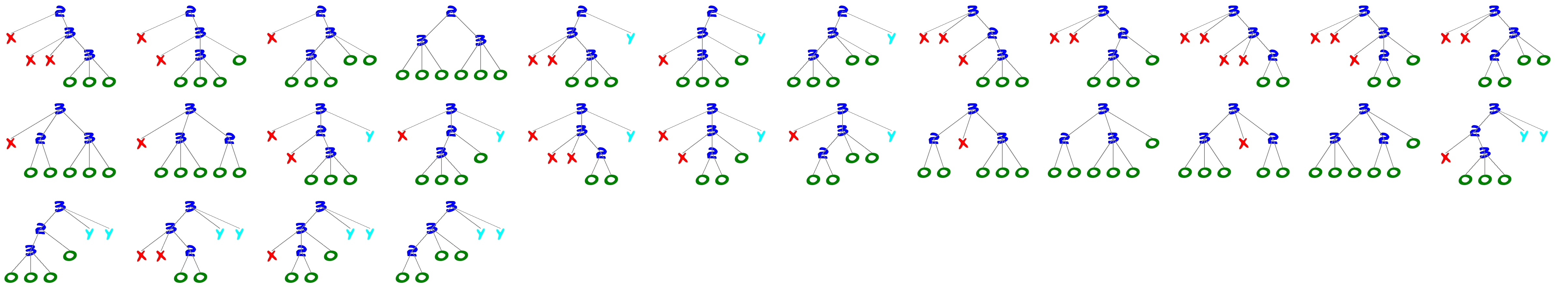} \caption{$C[1,2]=28, G[1,2]=110$}
\label{fig:DCG12}
\end{figure}
\begin{figure}[H]
\centering
    \includegraphics[width=1\linewidth]{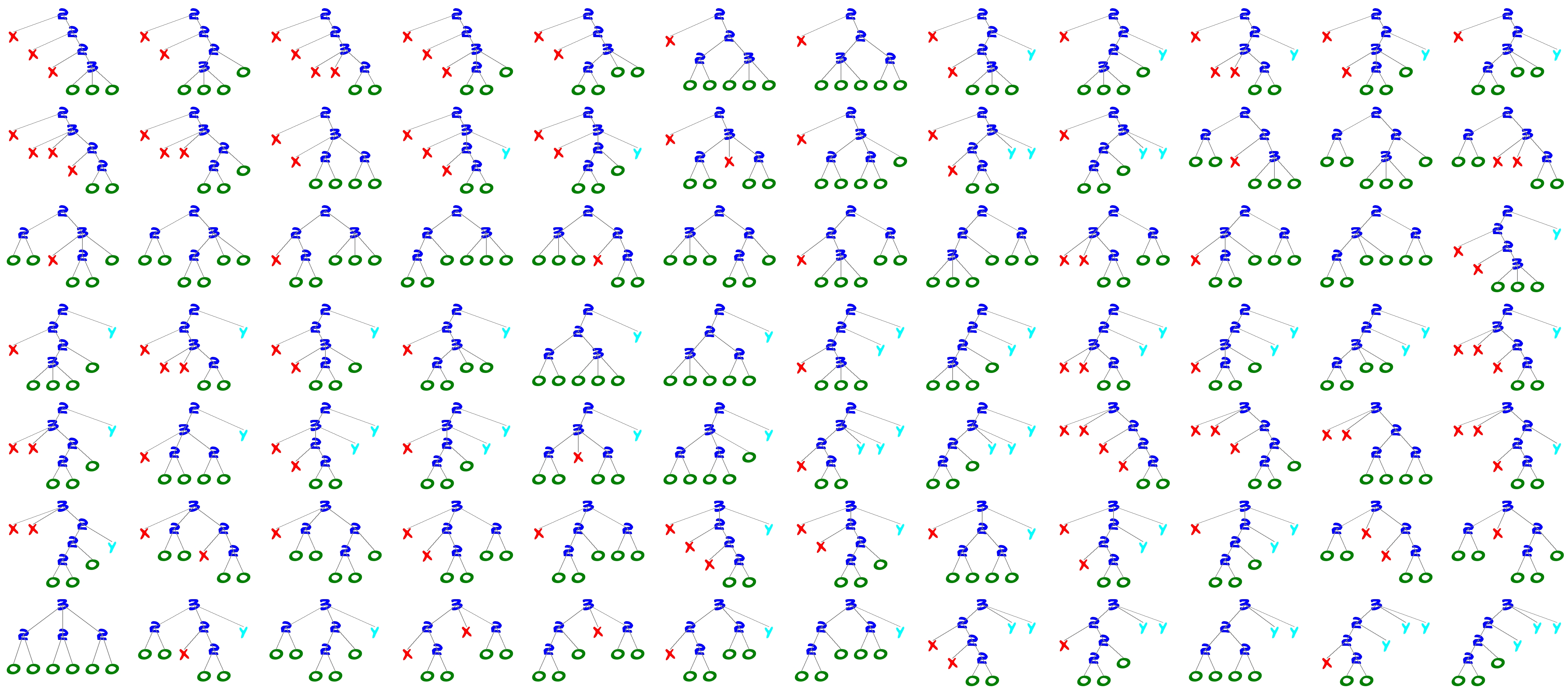} \caption{$C[3,1]=84, G[3,1]=288$ but conj. G3 counts 316}
\label{fig:DCG31}
\end{figure}

We failed hard on $G[3,1]$; we need a new conjecture.
\end{solution}
\begin{exercise}\label{ex:code}
Make tables of hyper-Catalan and Geode numbers.    
\begin{figure}[H]
\centering
    \includegraphics[width=1\linewidth]{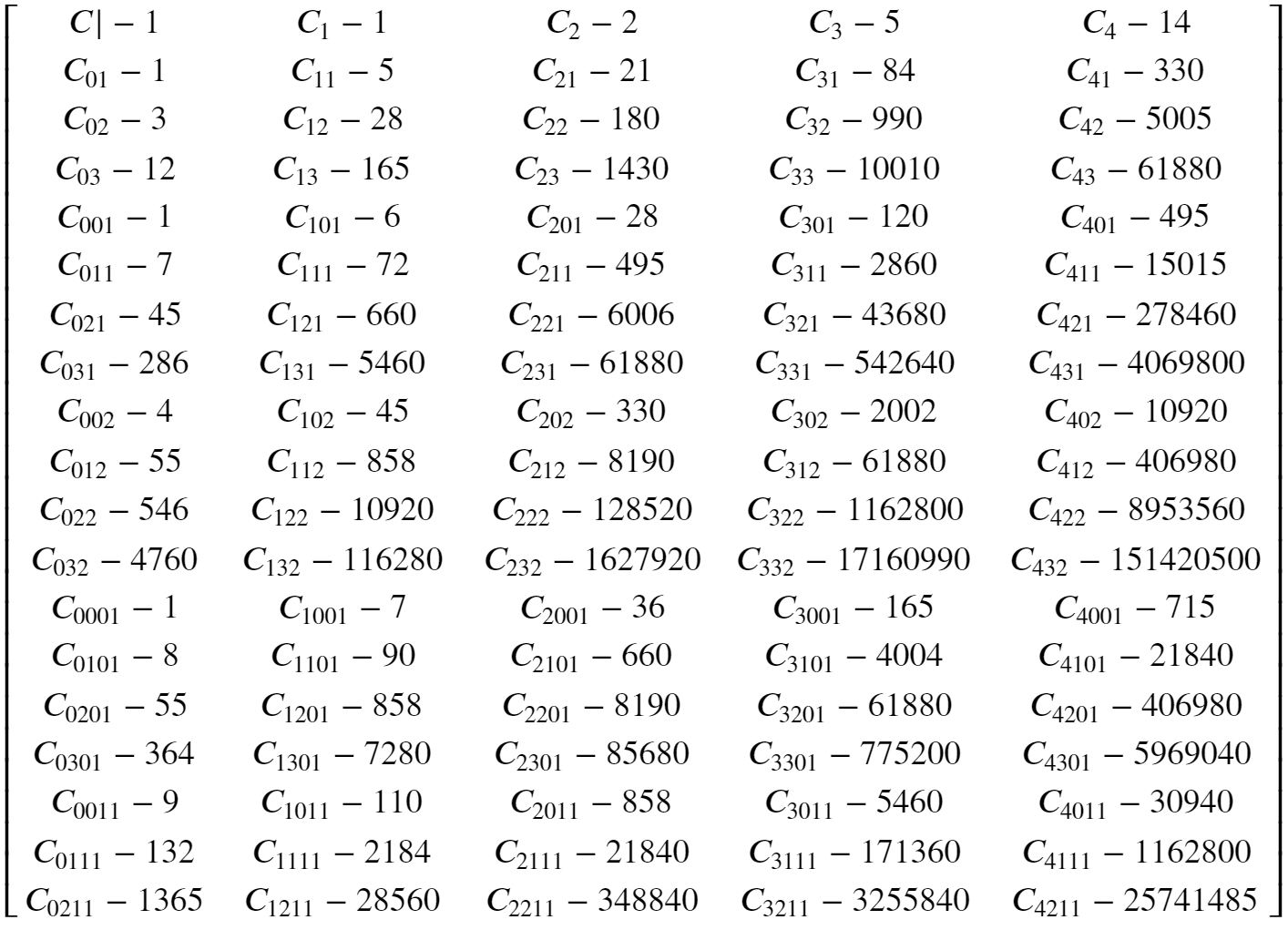} \caption{hyper-Catalan array}
\label{fig:HA}
\end{figure}
\begin{figure}[H]
\centering
    \includegraphics[width=1\linewidth]{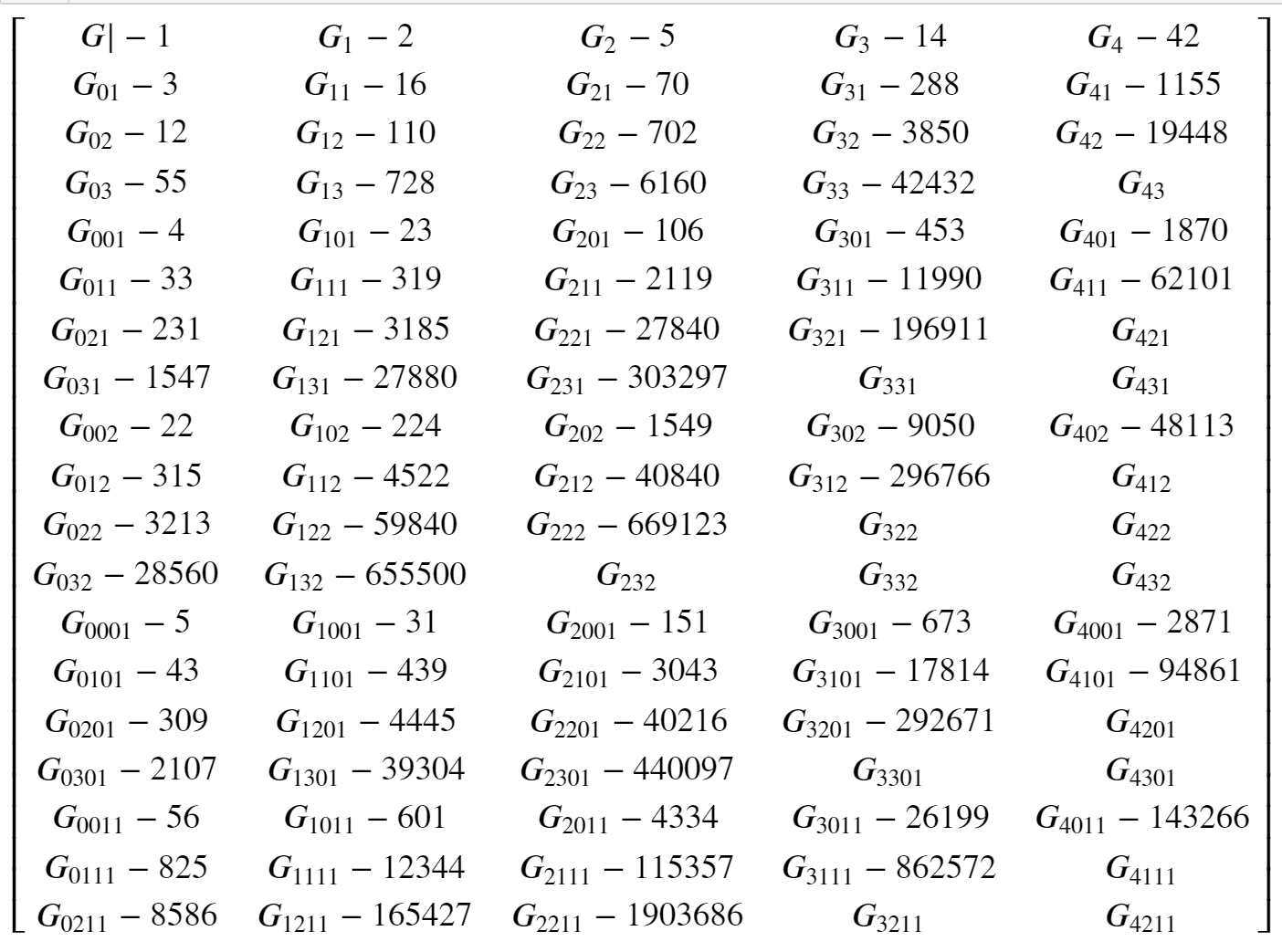} \caption{Geode array}
\label{fig:GA}
\end{figure}
\end{exercise}

\begin{solution}  
See the figures. 
We need slightly larger polynomials to pick up those missing Geode entries. 
Instead, let's work up some Python code to calculate Geode entries.
Let's start with some code to calculate the hyper-Catalan number $\mC$; the type vector $\m=[m_2, m_3, \ldots]$.
\begin{verbatim}
from sympy import *
def C(m):  # m=[m2,m3,...] so m2=m[0] etc.
    return factorial(E(m)-1) / (Fact(m) * factorial(V(m)-1))
def V(m):
    return 2 + sum((i+1)*m[i] for i in range(len(m)))
def E(m):
    return 1 + sum((i+2)*m[i] for i in range(len(m)))
def Fact(m):
    return prod(factorial(m[i]) for i in range(len(m)))
\end{verbatim}
For the Geode, we use a recurrence that turns each Geode entry into an integer combination of hyper-Catalans.  
See Rubine~\cite{Rubine2025} for details.
\begin{verbatim}
def G(m):
    if(m==[]):
        m = [0]
    maxi = list.index(m, max(m))
    m[maxi] += 1
    s = C(m)
    for i in range(0,len(m)):
        if i != maxi and m[i] > 0:
            m[i] -= 1
            s += -G(m)
            m[i] += 1
    m[maxi] -= 1
    return s
\end{verbatim}
\end{solution}

Let's test it out: $ G_{222}=669123 \quad\checkmark$.

OK, let's fill in the missing entries:
$ G_{232}=8754130,
 G_{331}=2580200, \newline 
 G_{322}=5826660,
 G_{332}=88952776,
 G_{3301}=3836840,
 G_{3211}=16990113, \newline 
 G_{43}=259350,
 G_{421}=1230250,
 G_{431}=18907196,
 G_{412}=1896650, \newline 
 G_{422}=43560528,
 G_{432}=762192600,
 G_{4201}=1872850,
 G_{4301}=28713476, \newline 
 G_{4111}=5648272,
 G_{4211}=129754776$.

\begin{exercise}
In June, 2024 someone sent us a YouTube video by Sateesh Mane (S.R. Mane),
who seems to know everything about the hyper-Catalan numbers, which he refers to as \textbf{Multi-parameter Fuss-Catalan Numbers}.
The corresponding paper is:
\textit{Multiparameter Fuss–Catalan numbers with application to algebraic equations} \cite{Mane2016}.

Mane first generalizes the two-parameter Fuss-Catalans to three-parameter Fuss-Catalan numbers, where the 1 becomes $r$, the power to which the solution series is raised.
By viewing the hyper-Catalans as a vector generalization of the three-parameter Fuss-Catalan numbers, 
Mane appears to have generalized our solution to equations with fractional exponents, and has an expression for the series solution raised to an arbitrary power.

The exercise is to transcribe Mane's solution from his paper 
and derive our solution from it.
\end{exercise}

\begin{solution} 
Mane writes the usual, two-parameter Fuss-Catalan numbers as
$$
A_t(m) = \dfrac{1}{(m-1)t+1}\binom{mt}{t} .
$$
\end{solution}
He generalizes this into the three-parameter Fuss Catalans,
$$
A_t(\mu, r) = \dfrac{r}{t!} \prod_{j=1}^{t-1} (t \mu + r - j) .
$$
which is a form he prefers as he's interested in the case where $\mu$ and $r$ are real (and complex).  
He notes the univariate generating function $B_\mu(r;z)$
$$
B_\mu(r;z) \equiv \sum_{t \ge 0} A_t(\mu, r) z^t
$$
satisfies 
$$
B_\mu(r;z) = (B_\mu(1;z))^r 
$$
and that $f=B_\mu(1;z)$
 is the solution to:
  $$f = 1 + zf^{\mu}.$$

Then he generalizes to solve:
\begin{equation} \label{eqn:mane}
f = 1 + z_1 f^{\mu_1} + \ldots + z_kf^{\mu_k}.
\end{equation}

For this he employs something hyper-Catalan-like, which he calls \textit{the multi-parameter Fuss-Catalan numbers}.
He uses a script $\mathcal{A}$; I won't bother with that, indicating the multi-parameter form with a bold vector subscript and a bold vector parameter.
\newcommand{\mub}{\boldsymbol\mu}
\newcommand{\zb}{\mathbf{z}}
$$
A_{\T}(\mub, r) = \dfrac{r}{t_1! \cdots t_k!} \prod_{j=1}^{-1 + \sum_i t_i} (\T \cdot \mub +r - j)
$$
We're using Mane's notation so his $\T$ is like our $\m$ and his $\zb$ is like our $\T$.
He defines the generating function:
$$
B(\mub; r; \zb) \equiv \sum_{\T \ge 0} A_{\T}(\mub, r) z_1^{t_1} \cdots z_k^{t_k}.
$$
His sum is over $\T \in \mathbf{N}^k$; I substituted our $\T \ge 0$ for that.
For convergence purposes he further specifies the order of summation in levels, what we would call face layers; I won't write that out.

He solves the above polynomial with real exponents and coefficients, equation (\ref{eqn:mane}), as $f=B(\mub; 1; \zb)$, and furthermore, $f^r=B(\mub; r; \zb)$.

Let's massage this into our solution.
We focus on the multi-parameter Fuss-Catalan numbers, and see how they relate to hyper-Catalans.
We're obviously interested in the $r=1$ case:
 $$
A_{\T}(\mub, 1) = \dfrac{1}{t_1! \cdots t_k!} \prod_{j=1}^{-1 + \sum_i t_i} (\T \cdot \mub +1 - j)
$$
Instead of the general exponents, let's solve the geometric polynomial, so exponents:
$$
\mub = [ 2, 3, 4, \ldots ].
$$
and corresponding coefficients
$$ \T = [t_1, t_2, \ldots] = [m_2, m_3, \ldots]=\m$$
so $t_i= m_{i+1}$.

$$
A_{ [m_2, m_3, \ldots]}([ 2, 3, 4, \ldots ], 1)  = \dfrac{1}{\m !} \prod_{j=1}^{-1 + \sum_{i \ge 2} m_i} (1-j+ 2m_2+3m_3 + \ldots)
$$
We're not concerned with the non-natural cases, so we can write the product as the ratio of factorials:
$$
A_{ [m_2, m_3, \ldots]}([ 2, 3, 4, \ldots ], 1)  = \dfrac{1}{\m !}  \dfrac{(2m_2+3m_3 + \ldots)!}{(1  - \left( \sum_{i \ge 2} m_i \right) + 2m_2+3m_3 + \ldots)!}
$$
$$
A_{ [m_2, m_3, \ldots] } ([ 2, 3, 4, \ldots ], 1)  = \dfrac{(2m_2+3m_3 + \ldots)!}{\m ! \, (1  + m_2+2m_3+3m_4 + \ldots)!} 
$$

Recall
$\mV = 2+\sum_{k \ge 2}  (k-1)m_k$,
$\ \mE= 1 + \sum_{k \ge 2}^{\fub} k \, m_k$,
$\ \mF = \sum_{k\ge 2} m_k$
and $
\mC =  \dfrac{(\mE-1)!}{(\mV-1) ! \, \m !} $.
We conclude:
$$
\mC = A_{ \m} ([ 2, 3, 4, \ldots ], 1)  
$$
so from Mane we can derive $f=\sum_{\m \ge 0}\mC \T^{\m}$ solves the geometric polynomial:
\begin{equation} \label{eqn:ggp}
f = 1 + t_2 f^{2} + \ldots + t_k f^k.
\end{equation}

We should be able to write powers $f^r$ using Mane's formulation.
We have $f= \sum_{\m \ge 0}\mC \T^{\m} = B([2,3,4,\ldots]; 1; \T)$ so $f^r =  B([2,3,4,\ldots]; r, \T)$.  
That suggests we generalize $ C_{\m}^{(r)} =   A_{\m}([2,3,4,\ldots], r)$  so $\mC= C_{\m}^{(1)}$.

\begin{align*}
 A_{\m}([2,3,4,\ldots], r) &= \dfrac{r}{m_2! \ m_3! \cdots } \prod_{j=1}^{-1 + \sum_i m_i} (r - j +  2m_2+3m_3 +4m_4 + \ldots)
 \\&
= \dfrac{r (r - 1 +  2m_2+3m_3 +4m_4 + \ldots) ! }{\m! (r+  m_2+2m_3 +3m_4 + \ldots)! } 
\\
  C_{\m}^{(r)}
&
= \dfrac{r (r-2 + \mE) ! }{\m! (r - 2 + \mV)! }  
\end{align*}
and the claim is:
$$
\left(\sum_{\m \ge 0}\mC \T^{\m} \right)^r = \sum_{\m \ge 0}  C_{\m}^{(r)} \, \T^{\m}.
$$

Cool, let's try some examples.
Start with a quadratic, $t=t_2, t_k=0$ for $k \ge 3$. 
For $r=2$ we already know how this one works; it's the Catalan numbers. We have $\m=[m]$, i.e. $m_2=m$.
The claim is:
$$
\left(\sum_{m \ge 0} C_m t^m   \right)^2
=\sum_{m \ge 0} C_{[m]}^{(2)} \, t^m.
$$
We know $\Tb =\sum_{m \ge 0} C_m t^m $ satisfies $\Tb = 1 + t\Tb^2$, i.e. $\Tb^2 = t^{-1}(\Tb -1)$ so we expect $C_{[m]}^{(2)}$ to be the shifted Catalan numbers, missing the initial one but otherwise the same.

$$
 C_{[m]}^{(2)} = \dfrac{2 (1 + 2m)!}{m!(2+m)!} = \dfrac{2(m+1) (-1 + 2(m+1))!}{(m+1)!(1+(m+1))!} 
$$
$$
 C_{[m]}^{(2)} = \dfrac{(2(m+1))!}{(m+1)!(1+(m+1))!} =C_{m+1} \quad\checkmark
$$

\begin{exercise} (From Irvin Miller, 20250614)

Solve $8x^3-6x+1=0$.

Irvin informs us:
$8x^3-6x+1= 8(x - \sin 10^\circ)(x - \sin 50^\circ)(x-\sin(-70^\circ)).$
\end{exercise}
\begin{solution}
Let's start out by verifying Irvin's claim.  Cosines are generally easier to work with; we seek a polynomial with zeros $\cos 80^\circ$, $\cos 40^\circ$, $\cos 160^\circ$, i.e. 1/9, 2/9 and 4/9 of a circle.

These angles satisfy $\cos 5 \theta = \cos 4 \theta$.  Proof:  That's $5 \theta = \pm 4 \theta + 360^\circ k $; minus sign solutions subsume the plus sign ones; we get $\theta = 40^\circ k$ for integer $k \quad \checkmark$.

So for $x=\cos \theta$ the equation $\cos 5 \theta = \cos 4 \theta$ becomes $T_5(x)=T_4(x)$ where those are the Chebyshev polynomials of the first kind, the multiple angle formulas for cosine.  That's 5th degree:
$16 x^5 - 20 x^3 + 5 x =8 x^4 - 8 x^2 + 1 $ or $16x^5-8x^4-20x^3+8x^2+5x-1=0$.  We know $x=\cos 0^\circ=1$ and $x=\cos 120^\circ=-1/2$ are solutions, so we factor out the corresponding factors giving:
$$
\dfrac{16x^5-8x^4-20x^3+8x^2+5x-1}{2(x- -1/2)(x - 1)}=8 x^3 - 6 x + 1 \quad\checkmark
$$

OK, Irvin has correctly identified the solutions.  Let's get some numbers out of our formula.  

It's a cubic; we proceed like in the paper by taking the generating series of the first four cross diagonals of the Bi-Tri array; we'll use the Geode factorization for fun:
$$
Q\left(  t_{2},t_{3}\right) =1+\left(  t_{2}+t_{3}\right) \left(
1+2t_{2}+3t_{3}+ 5t_{2}^{2}+16t_{2}
t_{3}+12t_{3}^{2}\right)
$$
To approximate $c_0 - c_1 x + c_2 x^2 + c_3 x^3=0$ we apply
$$
K\left(  c_{0},c_{1},c_{2},c_{3}\right)  \equiv\frac{c_{0}}{c_{1}} Q\left(  \frac{c_{0}c_{2}}{c_{1}^{2}},\frac{c_{0}^{2}c_{3}}{c_{1}^{3}}\right).
$$

We have $f(x)=1 - 6x + 8 x^3=0$ with approximate rational solutions:
\begin{align*}
\sin 10^\circ & \approx 0.17364817766693033 
\\[-4pt] \sin 50^\circ & \approx 0.766044443118978 
\\[-4pt] \sin(-70^\circ) & \approx -0.9396926207859083 \end{align*}

To apply our formula, we have $c_0=1, c_1=6, c_2=0, c_3=8$.  We get $K(1,6,0,8)=0.1736269877559315$ to sixteen decimal places, agreeing with $\sin 10^\circ $ to four decimal places.  To bootstrap we try:
$$f(0.1736+y)= 0.000254210047999921  - 5.27671296 y + 4.1664 y^{2}+8 y^{3}$$
and $K(0.0002542\ldots, 5.27671296, 4.1664, 8)=0.00004817766693033387$ so we get $x = 0.1736+y \approx 0.17364817766693033387$ which agrees to lots of decimal places.

What about the other roots?  We have to nudge the cubic we're solving toward the solution we want with an approximation.  Let's try $x\approx -1$ to aim at $\sin(-70^\circ)$.
$$f(-1+y) =- 1 - - 18 y - 24 y^{2}  + 8 y^{3} 
 $$
and $K(-1,-18,-24,8)=0.06028066015200691$ and when we add $-1$ we get $x\approx -0.939719339847993$, good to three decimal places. Taking $f( -0.939 + y)$ and applying $K$ to the coefficients gives $x\approx -0.9396926207858974$, again lots of agreement.

I'll leave $\sin 50^\circ$ for you.

\end{solution}

\begin{exercise}
Mats Granvik sent an email that Norman forwarded to me on 30 June 2025.
Mats included some Mathematica code that I couldn't unravel in the few minutes I looked at it.
But it's clear enough what's happening.
Mats is reverting series of the form:
\begin{equation*}
P_n(x) =  \dfrac{x}{\sum_{k=0}^n \binom{n}{k} x^k}   
\end{equation*}
and getting series with Fuss numbers as coefficients.  Can we show this?
\end{exercise}
\begin{solution}
We don't explicitly need binomial coefficients.  Letting $y=P_n(x)$, clearly:
$$ y = P_n(x)= \dfrac{x}{(1+x)^n}   $$
Let's substitute $z=1+x$ so $x=z-1$,
\begin{align*}
y &= \dfrac{z-1}{z^n}
\\ 1 - z + y z^n &= 0
\end{align*}
That's a trinomial that I believe first appeared in Lambert, 1758~\cite{Lambert1758}.  
For us, it's the geometric polynomial with a single non-zero coefficient, $y=t_n$.
Our general geometric polynomial solution is of course:
$$ 
0=1 - \alpha + t_2  \alpha ^2 + t_3  \alpha^3 + t_4 \alpha^4 + t_5 \alpha^5 \ + \  \ldots 
$$
has a formal power series solution
\begin{align*}
\alpha &= \Sb[t_2, t_3, \ldots] \equiv \!\!\!\! \sum_{m_2, m_3, \ldots \ge 0}  
\dfrac{( 2m_2 + 3m_3 + 4m_4 + \ldots )! }{(1 + m_2 + 2m_3 +\ldots)!\, m_2! \, m_3! \cdots}  \,   t_2^{m_2 } t_3^{m_3}\cdots .
\end{align*}
The solution to our equation, with only a single non-zero $t_n=y$, is $\Sb$ applied to $n-2$ zeros then $y$. Abbreviating $m=m_n$:
$$
z=\Sb[0,0,...,y] = \sum_{m \ge 0}  
\dfrac{( nm)! }{(1 + (n-1)m)!\, m! }  \,   y^{m}
$$
which are indeed the two-parameter Fuss numbers, see Lambert or Fuss, 1795~\cite{Fuss1795}.
We did essentially the same thing at the end of the Eisenstein section.  Oh, we have to remember: 
$$x=z-1 = \sum_{m > 0}  
\dfrac{( nm)! }{(1 + (n-1)m)!\, m! }  \,   y^{m}$$
\end{solution}

\vspace{10pt}
That concludes the exercises.
Please contact DeanRubineMath@gmail.com if you have questions or comments.

\vspace{10pt}
\noindent
%
\hrule

\end{document}